\documentclass[12pt]{amsart}
\usepackage{amsmath,amsfonts,amsthm,amscd,amssymb,euscript,
graphics,psfrag}
\usepackage{pictex}

\newtheorem{theorem}{Theorem}

\newtheorem{lemma}{Lemma}

\newtheorem{prop}{Proposition}

\newcommand{\trace}{{\rm trace}}

\newcommand{\Var}{{\rm Var}}

\newcommand{\nod}{{\rm mod}}
\newcommand{\orb}{\mathcal{O}}

\newcommand{\df}{\mathcal{F}}
\newcommand{\fd}{\mathcal{F}}
\newcommand{\gra}{\mathcal{G}}
\newcommand{\fK}{\mathcal{K}}

\newcommand{\bs}{\backslash}
\newcommand{\supp}{{\rm supp}}
\newcommand{\beq}{\begin{equation}}
\newcommand{\eeq}{\end{equation}}
\newcommand{\cx}{{\bbf C}}
\newcommand{\hyp}{{\bbf H}}
\newcommand{\disc}{{\bbf D}}
\newcommand{\lph}{{\Lambda (q) \bs \hyp}}
\newcommand{\loneh}{{\Lambda (1) \bs \hyp}}

\newcommand{\zed}{{\mathbb Z}}
\newcommand{\natls}{{\mathbb N}}
\newcommand{\ratls}{{\mathbb Q}}
\newcommand{\reals}{{\mathbb R}}

\newcommand{\fp}{{\mathbb F_p}}

\newcommand{\RR}{{\mathbb R}}

\newcommand{\var}{{\rm var}}

\newcommand{\q}{\zed/q \zed}

\newcommand{\p}{\zed/p \zed}

\newcommand{\cusp}{{\rm cusp}}
\newcommand{\flare}{{\rm flare}}

\newcommand{\SL}{{\rm SL}}

\newcommand{\SO}{{\rm SO}}

\def\bbf{\mathbb}
\def\mtrx#1#2#3#4{\begin{pmatrix} #1 & #2  \\ #3 & #4\end{pmatrix}}

\newcommand{\mL}{\mathcal{L}}

\newcommand{\mM}{\mathcal{M}}

\newcommand{\mF}{\mathcal{F}}

\newcommand{\SLtwofp}{\SL_2(\fp)}

\newcommand{\SLtwoZ}{\SL_2(\zed)}

\newcommand{\ga}{\alpha}

 \newcommand{\gd}{\delta}

\newcommand{\vge}{\varepsilon}

\newcommand{\gl}{\lambda}

\newcommand{\gr}{\varrho}

\newcommand{\gt}{\tau}

\numberwithin{equation}{section} \numberwithin{theorem}{section}

\def\vp{\varphi}
\def\arrowk{^\to{\kern -6pt\topsmash k}}
\def\arrowK{^{^\to}{\kern -9pt\topsmash K}}
\def\arrowr{^{^\to{\kern-}6pt\topsmash r}}
\def\arrowvp{^\to{\kern-8pt\topsmash\vp}}
\def\arrowf{^{^\to}{\kern -8pt f}}
\def\arrowg{^{^\to}{\kern -8pt g}}
\def\arrowu{^{^\to}a{\kern-8pt u}}
\def\arrowt{^{^\to}{\kern -6pt t}}
\def\arrowe{^{^\to}{\kern -6pt e}}
\def\tk{\tilde{\kern 1 pt\topsmash k}}
\def\barm{\bar{\kern-.2pt\bar m}}
\def\barN{\bar{\kern-1pt\bar N}}
\def\barA{\, \bar{\kern-3pt \bar A}}

\begin{document}

\title[Generalization of Selberg's  Theorem and Sieve]{Generalization of Selberg's $3/16$ Theorem and Affine Sieve}

\author{Jean Bourgain}
\address{School of Mathematics, Institute for Advanced Study, Princeton, NJ 08540} \email{bourgain@math.ias.edu}

\author{Alex Gamburd}
\address{Department of
Mathematics,  University of California at Santa Cruz, 1156 High
Street, Santa Cruz, CA 95064}
 \email{agamburd@ucsc.edu}

\author{Peter Sarnak}
\address{School of Mathematics, Institute for Advanced Study, Princeton, NJ 08540
and Department of Mathematics, Princeton University}
\email{sarnak@math.princeton.edu}

\dedicatory{Dedicated to the memory of Atle Selberg}

\thanks{The first author was supported in part by the NSF.
The second author was supported in part by DARPA, NSF, and the
Sloan Foundation. The third author was supported in part by the
NSF}

\maketitle

\section{Introduction}
A celebrated theorem of Selberg \cite{Se65} states that for
congruence subgroups of $\SLtwoZ$ there are no exceptional
eigenvalues below $3/16$.  We prove a  generalization of Selberg's
theorem for infinite index ``congruence'' subgroups of $\SLtwoZ$.
Consequently we  obtain sharp upper bounds in the affine linear
sieve, where in contrast to \cite{BGS} we use an archimedean norm
to order the elements.

Let $\Lambda$ be a finitely generated non-elementary subgroup of
$\SL_2(\zed)$; let $X_{\Lambda} =\Lambda \bs \hyp$ be the
corresponding hyperbolic surface (which is of infinite volume if
$\Lambda$ is of infinite index in
 $\SL(2, \zed)$). Let $\delta(\Lambda)$ denote the  Hausdorff
dimension of the limit set of $\Lambda$.
 The generalization of Selberg's theorem splits
into two cases: $\delta(\Lambda)> \frac{1}{2}$ and $0<
\delta(\Lambda) \leq \frac{1}{2}$.

 In the case that $\delta(\Lambda)> \frac{1}{2}$ the spectrum of the Laplace-Beltrami operator on
$L^2(X_{\Lambda})$ consists of finite number of points in $[0,
\frac{1}{4})$ (see \cite{LP82}). We denote them by $$0 \le
\lambda_0(\Lambda) < \lambda_1 (\Lambda) \le \dots \le
\lambda_{\max}(\Lambda) < \frac{1}{4}.$$  The assumption that
$\delta(\Lambda) > \frac{1}{2}$ is equivalent to
$\lambda_0(\Lambda) < \frac{1}{4}$, and in this case $\delta
(1-\delta) =\lambda_0$ \cite{Pa}.

The following extension of Selberg's theorem is proved in section
\ref{sec:t4}.

\begin{theorem} \label{t4}
Let $\Lambda$ be a finitely generated subgroup of $\SL(2, \zed)$
with $\delta(\Lambda)
> \frac{1}{2}$.  For $q \ge 1$ let $\Lambda(q)$ be the
``congruence'' subgroup $\{x \in \Lambda :  x \equiv I \mod q\}$.
There is $\varepsilon = \varepsilon(\Lambda) > 0$ such that
$$\lambda_1(\Lambda(q)) \ge \lambda_0(\Lambda(q)) +\varepsilon, $$
for all square-free $q \ge 1$  (note that
$\lambda_0(\Lambda(q))=\lambda_0(\Lambda)$).
\end{theorem}

 In \cite{ag} an explicit and
stronger version of Theorem \ref{t4} is proven under the
assumption that $\delta(\Lambda)
> \frac{5}{6}$.
See \cite{S05} for the sharpest known bounds towards Selberg's
$\frac{1}{4}$ Conjecture as well as bounds towards the Ramanujan
Conjectures for more general groups.

Theorem \ref{t4} is a consequence of Theorem 1.2 in \cite{BGS} and
the following result, which is of independent interest.

\begin{theorem} \label{t0} Let $\Lambda = \langle S \rangle$ be a
finitely generated subgroup of $\SL(2, \reals)$ with
$\delta(\Lambda) > \frac{1}{2}$.  Let $\{N_i\}$ be a family of
finite index normal subgroups of $\Lambda$.  Then the following
are equivalent
\begin{enumerate} \item The Cayley graphs $\gra(\Lambda/ N_i, S)$ form
a family of expanders.   \label{item1} \item There is $\varepsilon
= \varepsilon(\Lambda)
> 0$ such that
$\lambda_1(\Lambda/ N_i) \ge \lambda_0(\Lambda/ N_i) +\varepsilon
$.\label{item2}
\end{enumerate}
\end{theorem}

The argument in section \ref{sec:t4} establishes that $\ref{item1}
\Rightarrow \ref{item2}$;  the implication $\ref{item2}
\Rightarrow \ref{item1}$ is proved using Fell's continuity of
induction in section 7 of \cite{ag}.   Theorem \ref{t0}
generalizes results of Brooks \cite{brooks86} and Burger
\cite{burger86, burger88} who proved it in the case of co-compact
$\Lambda$.

\medskip

Combining Theorem \ref{t4} with Lax-Phillips theory of asymptotic
distribution of lattice points \cite{LP82} we obtain the following
result, which is the crucial ingredient in the execution of the
affine linear sieve in the archimedean norm.

\begin{theorem} \label{t4a} Let $\Lambda$ be a finitely generated subgroup of $\SL(2, \zed)$
with $\delta(\Lambda)
> \frac{1}{2}$.  Assume that $q$ is square free and $(q, q_0) =1$,
where $q_0$ is provided by the strong approximation theorem
\cite{mvw}.  There is $\varepsilon_1 > 0$ depending on $\Lambda$
such that for any $g \in \SL_2(q)$ we have \beq \label{et4a}
\begin{split} &|\{\gamma \in \Lambda \, | \, \|\gamma\| \leq T \,
\text{and} \, \gamma \equiv  g \, \nod \, q \} | \\ & = \frac{
c_{\Lambda} T^{ 2 \delta}}{|\SL_2(q)|} + O (q^3 T^{2 \delta -
\varepsilon_1}).
\end{split} \eeq
\end{theorem}

We now turn to the discussion of the case $\delta(X) \leq
\frac{1}{2}$. In this case there  is no discrete $L^2$ spectrum
and its natural replacement is furnished by the resonances of $X$,
which are given as the poles of the meromorphic continuation of
the resolvent $R_X(s) = (\Delta_X - s(1-s))^{-1}$. By the result
of Patterson \cite{Pa} and Sullivan \cite{DS82} $R_X(s)$ is
analytic for $\Re(s) > \delta$; Mazzeo and Melrose \cite{mm}
proved that $R_{X}(s)$ has a meromorphic continuation to the
entire plane.  In \cite{Pat88} Patterson proved that $R_{X}(s)$
has a simple pole at $s=\delta$ and no further poles on the line
$\Re(s) = \delta$; his proof is based on ideas from ergodic theory
related to Ruelle zeta-function.  Using further development of
these ideas due to Dolgopyat \cite{dolg}, Naud \cite{naud05a} has
recently established that   $R_X(s)$  is holomorphic (with the
exception of simple pole at $s= \delta$)  for $\Re(s)
> \delta - \varepsilon$, with $\varepsilon$ depending on $X$.
The following result, giving a resonance-free region for
congruence resolvent,  is proved in section \ref{s7a}.

\begin{theorem} \label{t5a}
Let $\Lambda$ be a finitely generated subgroup of $\SL(2, \zed)$
with $\delta(\Lambda)  \leq \frac{1}{2}$.  For $q \ge 1$ square
free  let $\Lambda(q)$ be the ``congruence'' subgroup $\{x \in
\Lambda : x \equiv I \mod q\}$; let $X(q) = \Lambda(q) \bs \hyp$.
There is $\varepsilon = \varepsilon(\Lambda)
> 0$ such that $R_{X(q)}(s)$ is holomorphic( with the
exception of simple pole at $s= \delta$)  for $\Re(s)
> \delta - \varepsilon \min \left(1, \frac{1}{\log(1 + |\Im
s|)}\right)$.
\end{theorem}

When $\delta \leq \frac{1}{2}$ we cannot apply the expansion
property \cite{BGS} directly; instead, to prove theorem \ref{t5a}
we use a dynamical treatment and invoke a generalization of the
underlying result on measure convolution (``$L^2$-flattening
lemma''): see Lemmas \ref{l1} and
 \ref{l1q} in section \ref{s2}.  It is likely that by
combining our methods with the extension of
 Dolgopyat's result \cite{dolg} to vector-valued functions,
analyticity of $R_{X(q)}(s)$ can be established for $\Re(s)
> \delta - \varepsilon$ --- in complete analogy\footnote{The analogy between Theorem \ref{t4} and Theorem \ref{t5a}
becomes clearer when their assertions are expressed in terms of
the Selberg zeta function \cite{Se56}.  If $\Lambda$ is a finitely
generated subgroups of $\SL_2(\reals)$ the Selberg zeta function
$Z_{X}(s)$  associated with $X = \Lambda \bs \hyp$ is known to be
an entire function, whose non-trivial zeros are given by the
resonances and the finite point spectrum \cite{glz, patper}.
Consequently, Theorem \ref{t4} is equivalent to the assertion that
when $\delta(\Lambda) > \frac{1}{2}$ there is
$\varepsilon(\Lambda)>0$ such that $Z_{X(q)}(s)$ is analytic and
non-vanishing on the set $\{\Re(s) > \delta - \varepsilon \}$,
except at $s=\delta$ which is a simple zero, while Theorem
\ref{t5a}  is equivalent to the assertion that when
$\delta(\Lambda) \leq \frac{1}{2}$ there is
$\varepsilon(\Lambda)>0$ such that $Z_{X(q)}(s)$ is analytic and
non-vanishing on the set $\left\{\Re(s) > \delta - \varepsilon
\min \left(1, \frac{1}{\log(1 + |\Im s|)}\right) \right\}$, except
at $s=\delta$ which is a simple zero.}
 in  with Theorem \ref{t4}.

Using methods of Lalley \cite{l89} we obtain  the following
analogue of Theorem \ref{t4a}, which is sufficient for sieving
applications.
\begin{theorem} \label{t5b} Let $\Lambda$ be a finitely generated subgroup of $\SL(2, \zed)$
with $0< \delta(\Lambda) \leq  \frac{1}{2}$.  Assume that $q$ is
square free and $(q, q_0) =1$, where $q_0$ is provided by the
strong approximation theorem \cite{mvw}.  There is $\varepsilon_1
> 0$, $C > 0$ depending on $\Lambda$ such that for any $g \in
\SL_2(q)$ we have \beq \label{et5a} \begin{split} &|\{\gamma \in
\Lambda \, | \, \|\gamma\| \leq T \, \text{and} \, \gamma \equiv g
\, \nod \, q \} | \\ & = \frac{ c_{\Lambda} T^{ 2
\delta}}{|\SL_2(q)|}\left(1 + O\left(T^{-\frac{1}{\log \log
T}}\right)\right) + O \left(q^C T^{2 \delta -
\varepsilon_1}\right).  \end{split} \eeq
\end{theorem}

We turn to applications to affine linear sieve \cite{BGS}.
 Consider the standard action on the two by two integer matrices by
multiplication on the left, and take the orbit $\orb$ of $I$ (the
identity matrix) under $\Lambda$. Set $|x|=\left(\sum_{i, j} x_{i
j}^2\right)^{\frac{1}{2}}$ where $x
=\mtrx{x_{11}}{x_{12}}{x_{21}}{x_{22}}$.  Set $N_{\Lambda}(T)=|\{x
\in \Lambda \, : \, |x| \le T\}|$ and let $\delta(\Lambda)$ be the
Hausdorff dimension of the limit set of an orbit $\Lambda z
\subset \hyp \cup \{\infty\}\cup \RR$, where $\hyp$ is the
hyperbolic plane, $z \in \hyp$ and $\Lambda$ act by linear
fractional transformations.  By the results of Lax-Phillips
\cite{LP82} and Lalley \cite{l89} we have  that $N_{\Lambda}(T)
\sim c_{\Lambda} T^{2 \delta(\Lambda)}$, as $T \to \infty$. Let $f
\in \ratls[x_{ij}]$ be integral on $\orb$ and assume that it is
weakly primitive for $\orb$, that is $\gcd \, \{ f(x) \, : \, x
\in \orb\}$ is $1$.  If $f$ is not weakly primitive then
$\frac{1}{N} f$ is, where $N = \gcd f(\orb)$, and we can represent
any weakly primitive $f$ as $\frac{1}{N} g$ with $g \in
\zed[x_{ij}]$ and $N= \gcd(\orb)$.

The coordinate ring $\ratls[x_{ij}]/(\det(x_{ij})-1)$ is a unique
factorization domain \cite{sans} and we can factor  $f$ into
$t=t(f)$ irreducibles $f_1 f_2 \dots f_t$ in this ring.  Set
$$\pi_{\Lambda, f}(T)=|\{x \in \Lambda; |x| \le T,
f_j(x) \, \text{is prime for} \, j=1, \dots, t \}|.$$

For $f \in \zed[x_{ij}]$  weakly primitive with $t(f)$ irreducible
factors our conjectured asymptotics is of the form: \beq
\label{conas} \pi_{\Lambda, f}(T) \sim \frac{c(\Lambda,
f)N_{\Lambda}(T)}{(\log T)^{t(f)}},  \, \, \text{as} \,  \, T \to
\infty, \eeq where $c(\Lambda, f)$ can be expressed as a product
of local densities; see \cite{fuchs, sarmaa} for an example of
explicit computation and numerical experiments.   In section
\ref{sec:t2} we establish the following sharp upper bound for $
\pi_{\Lambda, f}(T)$.

\begin{theorem} \label{t2} Let $\Lambda$ be a subgroup of  $\SL(2,
\zed)$ which is Zariski dense in $\SL_2$ and let $f \in
\zed[x_{ij}]$ be weakly primitive with $t(f)$ irreducible factors.
 Then
\beq \label{upb} \pi_{\Lambda, f}(T) \ll
\frac{N_{\Lambda}(T)}{(\log T)^{t(f)}}. \eeq
\end{theorem}

We also obtain the following lower bound for the number of points
$x \in \Lambda$ for which $f$ has at most a fixed number of prime
factors.

\begin{theorem} \label{t22} Let $\Lambda$ be a subgroup of  $\SL(2,
\zed)$ which is Zariski dense in $\SL_2$ and let $f \in
\zed[x_{ij}]$ be weakly primitive with $t(f)$ irreducible factors.
 Then there is an $r<\infty$, which can be
given explicitly in terms of $\varepsilon(\Lambda)$
 in theorems \ref{t4} and \ref{t5a},   such that
\beq \label{lowerb} |\{x \in \Lambda; |x| \le T, \, \, \text{and}
\, \, f (x) \, \text{has at most $r$ prime factors}  \}| \gg
\frac{N_{\Lambda}(T)}{(\log T)^{t(f)}}. \eeq
\end{theorem}

\section{Generalization of Selberg's $3/16$ theorem when $\delta > 1/2$}\label{sec:t4}
 Being a subgroup of finite index in $\Lambda(1)$,
$\Lambda(q)$ has the same bottom of the spectrum,
$\lambda_0(\lph)=\lambda_0(\loneh)$. As in section 2 of \cite{ag},
we have that for $q$ large enough $ \Lambda(1)/\Lambda(q) \cong
\SL_2(\q)$. Let $ S= \{A_1, \dots, A_k \},$ and let $S_q$ be the
natural projection of $S$ modulo $q$.

Theorem 1.7 in \cite{BGS} implies that if $\Lambda = \langle S
\rangle$ is non-elementary, then $\gra_q = \gra(\SL_2(\q), S_q )$
is a family of expanders. Consider the space $H(q)$  of
vector-valued functions $F$  on $\fd(1)$, satisfying \beq F(\gamma
z) = R_q(\gamma) F(z), \eeq for $\gamma \in \Lambda(1)/\Lambda(q)
\cong \SL_2(\q)$,  where $R_{q}(\gamma)$ denotes the regular
representation of $\SL_2(\q)$; we denote by $\langle , \rangle$
the inner product on this space and by $\| \,  \|$ the associated
norm.   Denoting by $H_0(q)$ the subspace of functions in $H(q)$
orthogonal to $\varphi_0$, the eigenfunction corresponding to
$\lambda_0$, the assertion of Theorem \ref{t4} is equivalent to
existence of $c>0$ such that \beq \frac{\int_{\fd}\|\nabla F \|^2
d \mu} {\int_{\df}\|F\| ^{2} d\mu} \geq  \lambda_0 +c \eeq for all
$F \in H_0(q)$.

Applying Theorem 1.7 in \cite{BGS}   for each $z \in \fd(1)$
  implies that
 there is $\varepsilon> 0$, depending \emph{only} on $S$, such that
for all $F \in H_0(q)$, we have  \beq \label{e:exp} \|F(\gamma z)
-F(z)\| \geq \varepsilon \|F(z)\| \quad \text{for some}\, \,
\gamma \in S . \eeq

Let $f=\|F\|$, and decompose it as \beq f=a\varphi_0(z)  +b(z),
\eeq where \beq \label{e:bort} \int_{\df} \varphi_0(z)
\overline{b(z)} d \mu(z) =0 \eeq and \beq \int_{\fd}|f|^{2}  d \mu
=a^{2} + \int_{\fd}|b|^2 d \mu  =1.\eeq Write $$F(z) = (F_1(z),
\dots, F_{N}(z)), $$ where $N =|\SL_2(\q)|$. Since
$$\nabla \left(\sum_{j=1}^{N}|F_j(z)|^2\right)^{\frac{1}{2}}=
\begin{cases} \frac{\sum_{j=1}^{N} F_{j}(z) \overline{\nabla F_{j}(z)}}{\left(\sum_{j=1}^{N}|F_j(z)|^{2}\right)^{\frac{1}{2}}} \, \, & \text{if} \, \,
\sum_{j=1}^{N}|F_j(z)|^2 \ne 0\\
0 \, \, & \text{otherwise} ,
\end{cases}
$$
we have
$$\|\nabla F\|^{2}(z)  \ge |\nabla \|F\||^2(z) = |\nabla f|^{2}(z).$$
Consequently we obtain: \beq
\begin{split} &\frac{\int_{\fd}\|\nabla F\|^2 d \mu} {\int_{\df}\|F\|
^{2} d\mu} \ge \frac{\int_{\fd}|\nabla \|F\| |^2 d \mu}
{\int_{\df}\|F\| ^{2} d\mu} =
\frac{\int_{\fd}|\nabla f |^2 d \mu} {\int_{\df}|f| ^{2} d\mu} =\\
& \frac{\int_{\fd}\langle \Delta f, f \rangle d \mu}
{\int_{\df}|f| ^{2} d\mu} = \int_{\fd}\langle a \lambda_0
\varphi_0 + \Delta b,   a \varphi_0 +b \rangle d \mu =\\& a^2
\lambda_0 + \langle \Delta b, b\rangle
\overset{\eqref{e:bort}}{\ge} a^2 \lambda_0 + \lambda_1
\int_{\fd}|b|^2 d \mu \ge \lambda_0 +(\lambda_1 -\lambda_0)
\int_{\fd}|b|^2 d \mu.
\end{split}
\eeq

By a theorem of Lax and Phillips \cite{LP82} there are only
finitely many discrete eigenvalues of $\Lambda$ in $[0,
\frac{1}{4}]$; consequently \beq \lambda_1 -\lambda_0 \geq c_1 > 0
.\eeq Therefore, as soon as $\int_{\fd}|b|^2 d \mu > \varepsilon_1
> 0$, we have that \beq \frac{\int_{\fd}|\nabla F |^2 d \mu}
{\int_{\df}|F| ^{2} d\mu} \geq \lambda_0 +c_1 \varepsilon_1. \eeq
Now consider the case of $\int_{\fd}|b|^2 d \mu =0$. We can assume
that  $a=1$ and write $F(z) =u(z) \varphi_0(z)$, with $u(z)
=(u_1(z), \dots, u_N(z))$, where $N =|\SL_2(\q)|$. Now \beq
\|u(z)\| = \sum_{j=1}^{N} |u_j|^2(z)= 1 \eeq implies \beq
\sum_{j=1}^{N} u_j \frac{\partial u_j}{\partial x} =
\sum_{j=1}^{N} u_j \frac{\partial u_j}{\partial x} =0, \eeq and
since
$$\frac{\partial (\varphi_0 u_j)}{\partial x} = u_j \frac{\partial
\varphi_0}{\partial x}  + \varphi_0 \frac{\partial u_j}{\partial
x}, $$
$$\frac{\partial (\varphi_0 u_j)}{\partial y} = u_j \frac{\partial
\varphi_0}{\partial y}  + \varphi_0 \frac{\partial u_j}{\partial
y}, $$ we have that \beq  \begin{split} & \|\nabla \varphi_0
u\|^2=\\
& \left(\frac{\partial \varphi_0}{\partial x}\right)^2
\sum_{j=1}^{N} u_j^2 + \varphi_0^2\sum_{j=1}^{N}\left(
\frac{\partial u_j}{\partial x}\right)^2 + \\
& \left(\frac{\partial \varphi_0}{\partial y}\right)^2
\sum_{j=1}^{N} u_j^2 + \varphi_0^2\sum_{j=1}^{N}\left(
\frac{\partial u_j}{\partial y}\right)^2 + \\
& 2 \varphi_0 \frac{\partial \varphi_0}{\partial x} \sum_{j=1}^{N}
u_j \frac{\partial u_j}{\partial x}+\\
& 2 \varphi_0 \frac{\partial \varphi_0}{\partial y} \sum_{j=1}^{N}
u_j \frac{\partial u_j}{\partial y}=\\
& |\nabla \varphi_0|^2 + \varphi_0^2 \|\nabla u\|^2.
\end{split}
\eeq

Consequently, \beq \begin{split} & \frac{\int_{\fd}\|\nabla F \|^2
d \mu} {\int_{\df}\|F\| ^{2} d\mu}=\frac{\int_{\fd}|\nabla
\varphi_0|^2 + \varphi_0^2 \|\nabla u\|^2 d \mu}
{\int_{\df}|\varphi_0| ^{2} d\mu} =\\ & \frac{\int_{\fd}|\nabla
\varphi_0|^2  d \mu} {\int_{\df}|\varphi_0| ^{2} d\mu} +
\frac{\int_{\fd} \varphi_0^2 \|\nabla u\|^2 d \mu}
{\int_{\df}|\varphi_0| ^{2} d\mu} \geq \lambda_0 +
\frac{\int_{\fd} \varphi_0^2 \|\nabla u\|^2 d \mu}
{\int_{\df}|\varphi_0| ^{2} d\mu}.
\end{split} \eeq
Our aim now is to show that \beq \label{e:grad} \frac{\int_{\fd}
\varphi_0^2 \|\nabla u\|^2 d \mu} {\int_{\df}|\varphi_0| ^{2}
d\mu} \geq c_2 > 0. \eeq To that end we assume that   that \beq
\label{e:grad1} \frac{\int_{\fd} \varphi_0^2 \|\nabla u\|^2 d \mu}
{\int_{\df}|\varphi_0| ^{2} d\mu} < \kappa \eeq and will obtain a
contradiction for sufficiently small $\kappa$  ($\kappa_j$ below
are of the form $a_j \cdot  \kappa$ for suitable constants $a_j$).
 Consider the fundamental domain $\fd(1)=\Lambda(1) \bs \disc$. Its
boundary, $\partial \df(1)$, consists of finitely many geodesic
arcs $\{l_i\}$ splitting into pairs $l_j$, $l_j'$ in such a way
that there is $\gamma_j \in S$ so that $l_j =\gamma_j l_j'$;
$\gamma_j$ are distinct and generate $\Lambda(1)$. Further, we
have decomposition of the following form:

\[
\df(1)=\fK(1)\cup \bigcup_{i \in Cu(1)} \cusp_{i}\cup \bigcup_{j
\in Fl(1)} \flare_{j}
\]

where
\begin{itemize}
\item[(1)] $\fK(1)$ is relatively compact in $\disc$

\item[(2)]$Cu(1)$ is a set of cusps of $\df(1)$. Each $\cusp_{i}$
is isometric to a standard cuspidal fundamental domain $P(Y_{i})$
of the form
\[
P(Y) = \{z=x+i y \mid  0<x<1, y>Y \},
\]
based on a horocycle
\[
h_{Y} = \{x+iy \mid y=Y \}.
\]

\item[(3)] $Fl(1)$ is a set of flares of $F(1)$.
 Each $\flare_{j}(\ga)$ is isometric  to a standard hyperbolic
fundamental domain $F(\alpha)$ of the form
 \[
F(\ga) =\left\{z: 1<|z|<\exp (\beta); 0< \arg (z) < \ga \right\},
\]
where $\alpha < \frac{\pi}{2}$.
\end{itemize}

Since  $\varphi_0 \in L^{2}(\fd(1))$, we have that \beq
\int_{\fK}|\varphi_0| ^{2} d\mu \geq c_3 \int_{\fd}|\varphi_0|
^{2} d\mu \,  \, \text{for some} \, \, c_3 > 0, \eeq and therefore
\eqref{e:grad1} implies that \beq \label{e:grad2} \frac{\int_{\fK}
\varphi_0^2 \|\nabla u\|^2 d \mu} {\int_{\fK}|\varphi_0| ^{2}
d\mu} \leq \kappa_1.  \eeq

We recall the definition of Fermi coordinates.   Let $\eta$ be the
geodesic in the hyperbolic plane parameterized with the unit speed
in the form
\[
t \rightarrow \eta(t)\, \in \hyp^2 \quad t\in \reals .
\]
Then $\eta$ separates $\hyp^2$ into two half-planes: a left hand
side and a right hand side of $\eta$.  For each $p \in \hyp^2$ we
have the directed distance $\rho$ from $p$ to $\eta$.  There
exists a unique $t$ such that the perpendicular from $p$ to $\eta$
meets $\eta$ at $\eta(t)$. Now $(\rho\, , t)$ is a pair of Fermi
coordinates of $p$ with respect to $\eta$.  In these coordinates
the metric tensor is
\begin{equation}  \label{e:fermi}
d s^{2} = d \rho^{2} +\cosh^{2}\rho d t^{2}.
\end{equation}

Introduce Fermi coordinates based on the bounding geodesics $l_j$,
and use them to foliate $\fK$. By compactness, using
\eqref{e:grad2}, we can find $z \in \fK$ and $\delta > 0$ such
that \beq \int_{B(z, \delta)}|\varphi_0|^2 > c_4 > 0, \eeq and for
all $j=1, \dots, k$

\beq \label{e:grad3} \frac{\int_{T_{j}(\delta)} \varphi_0^2
\|\nabla u\|^2 d \mu} {\int_{T_{j}(\delta)}|\varphi_0| ^{2} d\mu}
< \kappa_2, \eeq where $T_j(\delta)$ is a tube lying in $\fK$ and
containing $B(z, \delta)$ along the perpendicular to $l_j$.

Each $T$ is of the form $[-\delta, \delta] \times [\rho_{1, j},
\rho_{2, j}]$ in the appropriate Fermi coordinates.  Rewriting
\eqref{e:grad3} in Fermi coordinates \eqref{e:fermi}, and using
the fact that
$$\left(\frac{\partial u_j}{\partial \rho}\right)^2 +
\left(\frac{\partial u_j}{\partial t}\right)^2  \ge
\left(\frac{\partial u_j}{\partial \rho}\right)^2,$$ by Fubini's
Theorem we obtain \beq \int_{T_{j}(\delta)} \varphi_0^2 \|\nabla
u\|^2 d \mu \ge 2 \delta \int_{\rho_{1, j}}^{\rho_{2,
j}}\varphi_0^2 \|u'(\rho)\|^2 \cosh(\rho) d \rho. \eeq

 Let $L$ denote
the maximal length of the tubes $T_j$. Using the fact that if
$|u(\rho_1) -u_(\rho_2)| \ge C$  and $\rho_1 -\rho_2 \le L$ then
$\int_{\rho_1}^{\rho_2} |u'(\rho)|^2 d \rho  > C^2/L$ (since
$$C^2 \le \left(\int_{\rho_1}^{\rho_2} u'(\rho) d \rho \right)^2
\le  \left(\int_{\rho_1}^{\rho_2} 1 d \rho \right)
\left(\int_{\rho_1}^{\rho_2}| u'(\rho)|^2 d \rho \right) ),$$ we
obtain that \eqref{e:grad3} implies that  for all $j=1, \dots, k$
we have \beq \label{e:grad4} \int_{B(z, \delta)}\varphi_0^2
\|u(\gamma_j z) - u(z)\|  d \mu(z) < \kappa_3 \int_{B(z,
\delta)}\varphi_0^2 d \mu(z). \eeq

On the other hand, since $F(z) =u(z) \varphi_0(z)$ and
$\varphi_0(\gamma z) =\varphi_0(z)$ for all $\gamma \in
\SL_2(\q)$, \eqref{e:exp} implies that
 there is $\varepsilon(S)> 0$ independent of
 $q$, such that
 \beq \label{e:exp1}   \|u(\gamma z)
-u(z)\|
>\varepsilon(S) \quad \text{for some}\, \, \gamma \in S . \eeq
Applying mean-value theorem, we see that \eqref{e:grad4} implies a
contradiction with \eqref{e:exp1} once $\kappa$ is small enough
depending on $\varepsilon(S)$; consequently we have proved the
validity of \eqref{e:grad} and the proof of Theorem \ref{t4} is
complete.

The adaption of the preceding argument to proving the implication
$\ref{item1} \Rightarrow \ref{item2}$ of theorem \ref{t0} is
straightforward, as is the generalization of this result to higher
dimensional hyperbolic spaces: the theorem of Lax and Phillips, of
which we made crucial use in the first part of the argument, holds
for geometrically finite subgroups of $\SO(n,1)$ with Hausdorff
dimension of the limit set greater than $n/2$; the second   part
of the argument proceeds as above by restricting to compact part
of the fundamental domain and foliating it using Fermi
coordinates.   In particular, by combining the $\hyp^3$ analogue
of Theorem \ref{t0} with \cite{varju} and Theorem 6.3 in
\cite{BGS} we obtain the following theorem which has applications
to integral Apollonian packings \cite{fuchs, ko, sarmaa}.

\begin{theorem} \label{tApol} Let $\Lambda$ be a geometrically-finite subgroup of $\SL_2(\zed[\sqrt{-1}])$ with $\delta(\Lambda)>1$  and such that the traces of elements
of $\Lambda$ generate the field $\ratls(\sqrt{-1})$.  There is
$\varepsilon = \varepsilon(\Lambda)> 0$ such that
$$\lambda_1(\Lambda(\mathcal{A})) \ge \lambda_0(\Lambda(\mathcal{A})) +\varepsilon$$
as $\mathcal{A}$ varies over squarefree ideals in
$\zed[\sqrt{-1}]$.
\end{theorem}

\section{Counting lattice points for $\delta > \frac{1}{2}$} \label{sec:count}

Recall that the Poincar\'{e} upper half-plane model is the
following subset of the complex plane $\cx$:
\[
\hyp^2 = \left \{ z=x+iy \in \cx \mid y>0\right \},
\]
with the hyperbolic metric
\begin{equation}  \label{e:b1}
ds^2=\frac{1}{y^{2}}(dx^{2} +d y^{2}).
\end{equation}
The distance function on $\hyp^2$ is explicitly given by
\begin{equation}  \label{e:b2}
\rho(z,w)= \log \frac{|z-\bar{w}| +|z-w|} {|z-\bar{w}| -|z-w|}.
\end{equation}
We will  use the following expression:
\begin{equation}  \label{e:b3}
\cosh \rho (z, w) = 1+2 u(z, w),
\end{equation}
where
\begin{equation}  \label{e:b4}
u(z, w)=\frac{|z-w|^{2}}{4 \Im z \Im w}.
\end{equation}

The ring  $M_{2}(\reals)$ of two by two real matrices is a vector
space with inner product given by
\[
\langle g, h \rangle = \trace(g h^{t}).
\]
One easily checks that $\|g\|=\langle g,g\rangle^{\frac{1}{2}}$ is
norm in $M_{2}(\reals)$ and that
\begin{equation} \label{e:b6}
\|g\|^{2} = a^2+b^2 +c^2 +d^2 \quad \text{for} \, g=
\begin{pmatrix}
a & b \\
c & d
\end{pmatrix}.
\end{equation}
By taking $z=w=i$ in (\ref{e:b4})  we obtain that
\begin{equation} \label{e:b7}
\|g\|^{2} = a^2+b^2 +c^2 +d^2=4u(gi, i) +2.
\end{equation}

Now the result of Lax and Phillips \cite{LP82} is as follows.  Let
\beq N_{\Lambda}(T; z, w) = \#\{\gamma \in \Lambda \, : \, \rho(z,
\gamma w) \le T \}. \eeq Suppose $\delta > \frac{1}{2}$ and write
$\lambda_j=\delta_j(1-\delta_j)$; $\delta_0 = \delta$. Denoting
the eigenfunctions corresponding to $\lambda_j$ by $\varphi_j$ we
have \beq \label{e258} |N(T; z, w) - \sum_{j} c_j \varphi_j(z)
\overline{\varphi_{j}(w)}e^{\delta_j T}| = O(T^{5/6} e^{T/2}).\eeq

Turning to congruence subgroups $\Lambda(q)$ we have, using the
methods of  \cite{LP82},  that \beq \label{e258a}
|N_{\Lambda(q)}(T; z, w) - \sum_{j} c_j \varphi_{j,q}(z)
\overline{\varphi_{j,q}(w)}e^{\delta_{j,q} T}| = O(q^3 T^{5/6}
e^{T/2}),\eeq where implied constant is independent of $q$.

The base eigenfunction for $\Lambda(q)$, normalized to have $L^2$
norm one, is given by  \beq \varphi_{0, q} = \frac{1}{|\SL_2(\zed/
q\zed)|} \varphi_{0,1}. \eeq

Combining Theorem \ref{t4} with \eqref{e258a}  and \eqref{e:b2},
\eqref{e:b7} we obtain
 $\varepsilon_1 > 0$ depending on
$\Lambda$ such that for any $g \in \SL_2(q)$ we have $$ |\{\gamma
\in \Lambda \, | \, \|\gamma\| \leq T \, \text{and} \, \gamma
\equiv g \, \nod \, q \} | =  \frac{c_{\Lambda} T^{ 2
\delta}}{|\SL_2(q)|} + O (q^3 T^{2 \delta - \varepsilon_1}),$$
establishing Theorem \ref{t4a}.

\section{Shifts and thermodynamic formalism}\label{s1}
When $\delta \leq \frac{1}{2}$ the $L^2$ spectral theory of Lax
and Phillips \cite{LP82} is not available and we use symbolic
dynamics approach, in particular the work of Lalley \cite{l89}. In
this section we review the key necessary notions and results
pertaining to shifts of finite type.

A shift of finite type is defined as follows.  Let $A$ be an
irreducible, aperiodic $l \times l$ matrix of zeroes and ones,
called the transition matrix.  Define $\Sigma$ to be the space of
all sequences taking values in the alphabet $\{1, 2, , \dots, l\}$
with transitions allowed by $A$, that is

$$\Sigma=\{x \in \prod_{n=0}^{\infty} \{1, \dots, l\} \,  :
\, A(x_n, x_{n+1})=1 \, \, \forall n \}.$$

The space $\Sigma$ is compact and metrizable in the product
topology. Define the forward shift $\sigma: \Sigma \to \Sigma$ by
$(\sigma x)_n = x_{n+1}$ for $n \ge 0$.

Let $C(\Sigma)$ be the space of continuous, complex-valued
functions on $\Sigma$.  For $f \in C(\Sigma)$ and $0 < \rho < 1$
define

$$ \var_{n} f = \sup \left \{ |f(x) -
f(y)| \,  \,  : \,  x_j= y_j \, \text{for} \, 0 \leq j\leq n
\right \};$$

$$|f|_{\rho}  = \sup_{n\ge 0} \var_{n}(f)/ \rho^n;$$

$$\mathcal{F}_{\rho}= \{f \in C(\Sigma) \, : \, |f|_{\rho} <
\infty \}$$

Elements of $\mF_{\rho}$ are called Holder continuous functions.
The space $\mF_{\rho}$, when endowed with the norm $\|
\cdot\|_{\rho} = | \cdot|_{\rho} + \| \cdot\|_{\infty}$ is a
Banach space.

For $f, g \in C(\Sigma)$ define the transfer operator $\mL_{f}g
\in C(\Sigma)$ by
$$\mL_{f}g(x) = \sum_{y : \sigma y =x} e^{f(y)} g(y). $$
For each $\rho \in (0,1)$ and $f \in \mF_{\rho}$, $\mL_{f}:
\mF_{\rho} \to \mF_{\rho}$ is a continuous linear operator; if $f$
is real-valued then $\mL_{f}$ is positive.

Denoting
$$S_n f = f+ f \cdot \sigma + \dots f \cdot \sigma^{n-1},$$
we have
$$(\mL_{f}^n g)(x) = \sum_{y : \sigma^n y =x} e^{ S_n f(y)} g(y). $$

The following result is due to Ruelle; a proof can be found in
\cite{PP} or \cite{ruelle}.

\begin{theorem} \label{tr}  Let  $f \in
\mF_{\rho}$ be a real-valued function.
\begin{enumerate}
\item There is a simple eigenvalue $\lambda_{f} > 0$ of $\mL_{f}:
\mF_{\rho} \to \mF_{\rho}$ with strictly positive eigenfunction
$h_{f}$. \label{tr1} \item   The rest of the spectrum of $\mL_{f}$
is contained in $\{z \in \cx \, :\, |z| \leq \lambda_{f} -
\varepsilon \}$ for some $\varepsilon > 0$. \label{tr2} \item
There is a Borel probability measure $\nu_{f}$ on $\Sigma$ such
that $\mL_f^{*} \nu_f = \lambda_{f} \nu_f$. \label{tr3} \item If
$h_f$ is normalized so that $\int h_f d \nu_f =1$ then for every
$g \in C(\Sigma)$
$$\lim_{n \to \infty} \| \lambda_f^{-n} \mL_{f}^n g - (\int g d
\nu_f) h_f \|_{\infty} = 0$$ \label{tr4} \item There exist
constants $C_1, \varepsilon_1$ such that for all $g \in
\mF_{\rho}$ and for all $n$
$$\| \lambda_f^{-n} \mL_{f}^n g - (\int g d
\nu_f) h_f \|_{\rho} \leq C_1 (1 - \varepsilon_1)^n \|
g\|_{\rho}.$$ \label{tr5}
\end{enumerate}
\end{theorem}

The pressure functional is defined by
$$P(f)= \sup_{\nu} \int f d \nu + H_{\nu}(\sigma),$$
where the supremum is taken over the set of $\sigma$-invariant
probability measures and $H_{\nu}(\sigma)$ is the measure
theoretic entropy of $\sigma$ with respect to $\nu$. We have  (see
\cite{PP} or \cite{ruelle})
$$P(f) = \log \lambda_f.$$

 A measure $\mu$ is called the
equilibrium state or the Gibbs measure with the potential $f$ if
$$\int f d \mu + H_{\mu}(\sigma) = P(f).$$
For $f \in \mF_{\rho}$  Gibbs measure $\mu_f$ is the unique
$\sigma$- invariant probability measure on $\Sigma$ for which
there exist constants $0 < C_1 \leq C_2 < \infty$ such that
$$C_1 \leq \frac{\mu_f \{ y \in \Sigma : y_i = x_i, 0 \leq i < n
\}}{ \lambda_f^{-n} \exp\{S_n f(x)\}} \leq C_2. $$

As will become clear in the next section, the analyticity
properties of the map $z \to \mL_{z f}$, $z \in \cx$ will play
crucial role in the proof. For $f \in \mF_{\rho}$ fixed,
real-valued function, such that $S_m f$ is strictly positive for
some $m$, the quantities $\mL_{zf}$, $\lambda_{zf}$, $h_{zf}$,
$\nu_{z f}$ will be abbreviated by $\mL_{z}$, $\lambda_{z}$,
$h_{z}$, $\nu_{z }$.

\section{Resolvent of transfer operator and  lattice count problem}

 Let $\Lambda$ be a Fuchsian group
with no parabolic elements (this condition is automatically
satisfied in the case $\delta(\Lambda) \leq  \frac{1}{2}$ -- see
\cite{LP82}), generated by $k$ elements  $g_1, \dots, g_k  \subset
\SL_2(\zed)$. We identify $\Lambda$ with $\Sigma_*$, defined as
the set of finite sequences in the alphabet $\{g_1, g_1^{-1},
\dots, g_k, g_k^{-1}\}$ $(l=2k)$ with admissible transitions.
According to Series \cite{series}, this may be done so as to
obtain a shift of finite type.

 Let  $w \in \hyp (=\disc)$, and suppose it is not a fixed point for any
$\gamma \in \Lambda$; let $d_{H}$ denote hyperbolic distance. For
$x= x_1, \dots, x_m \in \Sigma_* (=\Lambda)$ define \beq
\label{e51} \tau(x)= d_H(0, x_1 \dots x_m w)- d_H(0, x_2 \dots x_m
w). \eeq

The left shift $\sigma$ on $\Sigma$ corresponds to the Nielsen map
(see \cite{series})  $F : L \to L$, where $L$ denotes the limit
set of $\Lambda$.

Recalling that  \beq \label{e52} S_n \tau = \tau + \tau \cdot
\sigma + \dots + \tau \cdot \sigma^{n-1}, \eeq we have
 \beq \label{e53} S_n \tau (x) = d_H(0, x_1 \dots x_n x_{n+1}  \dots w) - d_H(0, x_{n+1} \dots w). \eeq

For $a \in \reals$, $x \in \Sigma_{*}$, $\phi: \Sigma_* \to
\reals$, let  \beq \label{e54} N(a, x) = \sum_{n=0}^{\infty}
\sum_{y: \sigma^n y =x} \phi(y) 1_{ \{S_n \tau (y) \leq a \} }.
\eeq

Clearly \beq \label{e55} N(a,x) = \sum_{y \in \Lambda, d_H(0, y x
w) - d_{H}(0, x w)  \leq a} \phi(yx), \eeq where in the summation
$y$ is restricted so as to make $yx$ admissible.

In particular, for $\phi=1$
$$N(a,x)= |\{ \gamma \in  \Lambda \,  | \, d_{H}(i, \gamma x w) -
d_{H}(i, x w) \leq a \} |. $$

Returning to \eqref{e54}, one has the renewal equation \beq
\label{e57} N(a,x) = \sum_{\sigma(x')=x} N(a-\tau(x'), x') +
\phi(x) 1_{\{a \geq 0\}} \eeq (cf. \cite[(2.2)]{l89}).

The link with the transfer operator comes by taking the Laplace
transform of \eqref{e57}. Defining for $\Re z < -C$ \beq
\label{e58} F(z,x) = \int_{-\infty}^{\infty}e^{az} N(a,x) d a ,
\eeq equation \eqref{e57} gives the relation \beq \label{e59}
F(z,x)= \sum_{\sigma(x')=x} e^{z \tau(x')} F(z, x') +
\frac{\phi(x)}{z}. \eeq Thus we have  \beq \label{e510} (I- \mL_z)
F(z,x) = \frac{\phi(x)}{z}. \eeq

This leads  us to the study of the resolvent $(I- \mL_z)^{-1}$.
Before stating   the results of Lalley \cite{l89} and Naud
\cite{naud05a} for $\mL_z|_{\mathcal{F}_{\rho}(\Sigma)}$ (Theorem
\ref{ln} below)  we recall the following reinterpretation of the
Hausdorff dimension of the limit set in terms of the pressure
functional (see \cite{PP} or \cite{ruelle}).   Because $\tau$ is
eventually positive, the variational principle implies (see
\cite{PP} or \cite{ruelle}) that the pressure functional $P(- x
\tau)$ is strictly decreasing and has unique positive zero. Define
$\delta$ by $P(- \delta \tau) =0$, that is, \beq \label{e41}
\lambda_{-\delta} =1. \eeq

\begin{theorem} \label{ln}
\begin{enumerate}
\item There is $\varepsilon > 0$ such that for $\Re z < - \delta +
\varepsilon$, $ z \notin U$ ($U$ a suitable neighborhood of $-
\delta$), we have \beq \label{e42}
\|\mL_z|_{\mathcal{F}_{\rho}(\Sigma)}\|_{\rho} \lesssim |\Im z|^2
e^{-\varepsilon n}. \eeq  \label{ln1} \item  For $ z \in U$
decompose on $\mathcal{F}_{\rho}(\Sigma)$ \beq \label{e43} \mL_z
=\lambda_z(\nu_z \otimes h_z) + \mL_{z}'', \eeq (where $z \to
\lambda_z, z \to h_z, z \to \nu_z$ are holomorphic extensions to
$U$, satisfying
$$\mL_z h_z= \lambda_z h_z, \quad  \mL^{*}_z \nu_z = \lambda_z \nu_z,
\quad \int h_z d \nu_z =1).$$ Then \beq \label{e44}
\|(\mL_{z}'')^n\|_{\rho} < e^{-\varepsilon n } \quad \text{for}
\quad z \in U. \eeq  \label{ln2}
\end{enumerate}
\end{theorem}

Part \ref{ln1} follows from the discussion in Lalley \cite[p.
25]{l89} (in the case when $\tau$ is non-lattice) and Theorem 2.3
in Naud \cite{naud05a} to provide \eqref{e42} when $|\Im z|$ is
large. Naud's work build crucially on the approach of Dolgopyat
\cite{dolg}.   Note that Lalley does not give explicit estimates
on $\| \mL_z^n\|$ for $\Re z = -\delta$ and $\Im z \to \infty$ and
certainly no bound of the strength of Theorem 2.3 in Naud
\cite{naud05a}.

Part \ref{ln2}  is Proposition 7.2. in \cite{l89}

\section{Lattice count in congruence subgroups for $\delta \leq \frac{1}{2}$}
In this section we modify the setup discussed in the preceding two
sections to the setting of congruence subgroups of $\Lambda$ ---
the modification is analogous to the one preformed in the case
$\delta
> \frac{1}{2}$.

Fix modulus $q$ such that $\pi_q(\Lambda)=\SL_2(q)$. Instead of
considering functions on $\Sigma$ (as in Lalley \cite{l89} ) we
consider functions on $\Sigma \times \SL_2(q)$.

For $f \in C(\Sigma \times \SL_2(q))$ define

$$ \|f\|_{\infty} = \max_{x} \left( \sum_{g \in \SL_2(q)}
|f(x,g)|^{2}\right)^{\frac{1}{2}};$$

$$ \Var_{n} f = \sup \left \{ \left( \sum_{g} |f(x,g) -
f(y,g)|^2\right)^{\frac{1}{2}} \, |\, x_j= y_j \, \text{for} \,
j\leq n \right \};$$

$$|f|_{\rho}=\sup_{n} \frac{\Var_{n}(f)}{\rho^n}.$$

Let $\mF_{\rho} = \mF_{\rho}(\Sigma \times \SL_2(q))$ denote the
space of $\rho$-Lipschitz continuous functions with the norm $$\|
\cdot\|_{\rho}=\|\cdot\|_{\infty}+|\cdot|_{\rho}.$$

Let $\tau: \Sigma_{*}\to \reals$ be given by \eqref{e51} and
consider the ``congruence transfer operator''  $\mM_z = \mM_{z
\tau}$ on $\mF_{\gr}(\Sigma \times \SL_2(q))$: \beq \label{e11}
\mM_{z} f(x, g)=\sum_{i=1}^{l} e^{z \tau(i,x)} f(i,x; g_i g),\eeq
where $z \in \cx$ and the summation is restricted so as to make
$(i, x)$ admissible.

Thus our $\mM_{z}$ differs from the one considered in \cite{l89}
in that it acts on functions on $\Sigma \times \SL_2(q)$ rather
than on functions on $\Sigma$:  the reason behind this difference
is the same as in the proof of the spectral gap when
$\delta(\Lambda)
> \frac{1}{2}$.

We have  \begin{equation*}\begin{split} \mM_{z}^2 f(x, g) &=
\sum_{i_1 =1}^{l} e^{z \tau(i_1, x)} (\mM_{z} f)((i_1, x), g_{i_1}
g)\\
&=\sum_{i_1, i_2 =1}^l e^{z \gt(i_1, x) +\gt(i_2 i_1, x))} f((i_2,
i_1,  x); g_{i_2} g_{i_1} g), \end{split} \end{equation*} and in
general for the $n$-th iterate we have \beq \label{e12} \mM_{z}^n
f(x; g) =\sum_{i_1, \dots, i_n =1}^l e^{z (\gt(i_n,\dots , i_1, x)
+ \gt(i_{n-1},\dots , i_1,  x) + \dots + \gt(i_1, x))} f((i_n
\dots i_1,  x); g_{i_n} \dots g_{i_1} g), \eeq where again the
summation is restricted to admissible words.

 From Ruelle's theorem
(Theorem \ref{tr}) it follows that \beq \label{e14} \sum_{i_1,
\dots, i_n =1}^l e^{ \Re z \left( \gt(i_n,\dots , i_1,  x) +
\gt(i_{n-1},\dots , i_1, x) + \dots + \gt(i_1, x)\right)} \sim
\lambda_{\Re z}^{n} \eeq for large $n$.

Let $\varphi$ be a function on $\SL_2(q)$.  Returning  to
 \eqref{e54}, \eqref{e55} we let
 \beq \label{e6.5}
 N(a,x) = \sum_{y \in \Lambda,
d_H(0, y  x w) - d_{H}(0, x w)  \leq a} \varphi(\pi_q(y)) \eeq and
$F(z,x)$ its Laplace transform defined by \eqref{e58}.

Then \beq \label{e6.6} F(z,x) = f(z, x, \pi_{q}(x)),\eeq where
$f(z, x, g)$ satisfies \beq \label{e6.7} (1- \mM_z)f = \frac{1
\otimes \varphi}{z} \eeq and $\mM_z$ is the congruence transfer
operator introduced above  (note that obviously $1 \otimes
\varphi$ is in $\mF_{\gr}(\Sigma \times \SL_2(q))$).

Our aim is to evaluate \beq \label{e56} N(a) =\sum_{y \in \Lambda,
d_H(0, y w) - d_{H}(0, w)  \leq a} \varphi(\pi_q(y)), \eeq which
gives the sum of $\varphi$ on the $\nod \, q$ reduction of the
hyperbolic ball
$$\{y \in \Lambda \, | \,  d_H(0, y w) - d_{H}(0, w)  \leq a\}.$$

Our goal now is to obtain the appropriate extension of Theorem
\ref{ln} to the setting of congruence subgroup.  As is to be
expected, it is at this point that expansion property will play a
crucial role.

\section{Expansion and $L^2$ flattening}\label{s2}   Let $\mu$ be a symmetric measure on $G=\SL_2(p)$ (for the sake of exposition we first consider the simpler case of prime $p$)  and
consider the convolution map $T: L^2(G) \to L^2(G)$, given by
$\varphi \to \mu \ast \varphi$.  Decomposing the regular
representation of $G$ into irreducible representations if follows
from the result of Frobenius \cite{Fr} that each eigenvalue
$\lambda$ of the convolution restricted to $L_0^2(G)$ occurs with
multiplicity at least $\frac{p-1}{2}$. Trace calculation yields
therefore $$|G| \|\mu\|_2^2= \sum_{x\in G} \langle T^2 \delta_x,
\delta_x \rangle \ge \frac{p-1}{2} \lambda^2.$$  Hence \beq
\label{e21} |\lambda| \leq \sqrt{\frac{2}{p-1}} \|\mu\|_2
|G|^{\frac{1}{2}}. \eeq

Recall also the $L^2$- flattening lemma proven in \cite{bg}. Let
$\mu \in \mathcal{P}(G)$ satisfy \beq \label{e22} \|\mu\|_{\infty}
< p^{-\tau} \eeq for some $\tau > 0$ and also \beq \label{e22a}
\mu(a G_1) < p^{-\tau} \eeq for all cosets of proper subgroups
$G_1$ of $G$. Given $\kappa
>0$ there is $l=l(\tau, x) \in \zed_{+}$ such that \beq
\label{e23} \|\mu^{(l)}\|_{2} < |G|^{-1/2 + \kappa}. \eeq

Denote $\mu'(x)=\mu(x^{-1})$.  Since $\mu \ast \mu'$ also
satisfies \eqref{e22},  \eqref{e22a} we have by \eqref{e23} that
\beq \label{e24} \|(\mu' \ast \mu)^{(l)}\|_{2} < |G|^{-1/2 +
\kappa}. \eeq Consider the convolution operator $T \varphi = \mu'
\ast \mu \ast \varphi$ and let $\lambda$ be an eigenvalue of $T$
on $L_0^2(G)$. Hence $\lambda^l$ is an eigenvalue of $T^l$ on
$L_0^2(G)$ and applying \eqref{e21} with $\mu$ replaced by $(\mu'
\ast \mu)^{(l)}$ implies
$$|\lambda|^l \leq \sqrt{\frac{2}{p-1}} \|(\mu' \ast \mu)^{(l)}\|_2
|G|^{\frac{1}{2}} < \sqrt{\frac{2}{p-1}} |G|^{\kappa} <
p^{-\frac{1}{4}}$$ if we take $\kappa < \frac{1}{4}$. Consequently
\beq \label{e25} |\lambda|< p^{-\frac{1}{4l}}. \eeq This means
that if $\varphi \in L_0^2(G)$, then
$$ \|\mu' \ast \mu \ast \varphi\|_2 \leq p^{-\frac{1}{4l}}
\|\varphi\|_2$$ and hence \beq \label{e26} \|\mu \ast \varphi\|_2
\leq p^{-\frac{1}{8l}} \|\varphi\|_2.\eeq We proved that if $\mu
\in \mathcal{P}(G)$ satisfies \eqref{e22}, \eqref{e22a} for some
$\tau
> 0$, then \beq \label{e27} \|\mu \ast \varphi\|_2 \leq p^{-\tau'}
\|\varphi\|_2, \quad \varphi \in L_0^2(G)\eeq for some $\tau' >0$.
Therefore if $\mu \in \mathcal{M}_{+}(G)$ satisfies \beq
\label{e28} \|\mu\|_{\infty} < p^{-\tau} \|\mu\|_1 \eeq and \beq
\label{e28a} |\mu|(a G_1) < p^{-\tau} \|\mu\|_1 \eeq for cosets of
proper subgroups $G_1$ then \beq \label{e29} \|\mu \ast
\varphi\|_{2} \leq p^{-\tau'} \|\mu\| \|\varphi\|_2, \quad \varphi
\in L_0^2(G). \eeq More generally, let $\mu \in \mathcal{M}(G)$
and decompose $\mu = \mu_{+} + \mu_{-}$. Estimate
$$\|\mu \ast \varphi\|_2 \leq \|\mu_+ \ast \varphi\|_2 + \|\mu_- \ast
\varphi\|_2.$$ Assume $\mu$ satisfies \eqref{e28} and
\eqref{e28a}. If $\|\mu_+\|> p^{-\frac{\tau}{2}} \|\mu\|$, we have
$$\|\mu_+\|_{\infty} \leq \|\mu\|_{\infty} <
p^{-\frac{\tau}{2}}\|\mu_+\|_{1} \quad \text{and} \quad \mu_{+}(a
G_1) < p^{-\frac{\tau}{2}} \|\mu_+\|_1$$ and \eqref{e29} implies
that for $\varphi \in L_0^2(G)$ we have
$$\|\mu_+ \ast \varphi\|_2 \leq p^{-\tau'}\|\mu_+\|\|\varphi\|_{2}
\leq p^{-\tau'}\|\mu\| \|\varphi\|_2. $$ If $\|\mu_+\| \leq
p^{-\frac{\tau}{2}}\|\mu\|$, then obviously
$$\|\mu_+ \ast \varphi\|_2 \leq \|\mu_+\|\|\varphi\|_{2}
\leq p^{-\frac{\tau}{2}} \|\mu\|\|\varphi\|_2. $$ Hence
\eqref{e29} holds.

The same considerations apply for $\mu \in \mathcal{M}_{\cx}(G)$;
consequently we obtain the following result.
\begin{lemma}\label{l1}  Given $\kappa>0$, there is $\kappa'>0$ such that if
$\mu \in \mathcal{M}_{\cx}(G)$ satisfies \beq \label{e210}
\|\mu\|_{\infty} < p^{-\kappa} \|\mu\|_{1} \eeq and \beq
\label{e210a} |\mu|(a G_1) < p^{-\kappa} \|\mu\|_1 \eeq for cosets
of proper subgroups $G_1$ of $G$, then \beq \label{e211} \|\mu
\ast \varphi\|_2 \leq C p^{-\kappa'} \|\mu\|_{1} \|\varphi\|_2
\quad \text{for} \quad \varphi \in L_0^2(G). \eeq Here $p$ is
assumed to be sufficiently large.
\end{lemma}

We have  a similar result for   $G=\SL_2(q)$ with $q$ square-free
(see \cite{BGS} and \cite{varju}). We make the following
decomposition of the space $L^2(\SL_2(q))$. For $q_1|   q$, define
$E_{q_1}$ as the subspace of functions defined $\mod q_1$ and
orthogonal to all functions defined $\mod q_2$ for some $q_2
|q_1$, $q_2 \neq q_1$. Hence \beq \label{e82} L^2(\SL_2(q))=\reals
\oplus \bigoplus_{q_1 |q} E_{q_1},\eeq which is, in fact, the
generalized Fourier-Walsh decomposition corresponding to the
product representation
$$\SL_2(q) \cong \prod_{p|q} \SL_2(p).$$
Let $P_1$ (respectively $P_{q_1}$) be the projection operator on
the constant functions (respectively $E_{q_1}$).

\begin{lemma}\label{l1q} Let $q$ be square free and $G=\SL_2(q)$.  For $\mu \in \mathcal{M}_{\cx}(G)$ and $q_1|q$  define
$|||\pi_{q_1}(\mu)|||_{\infty}$ to be the maximum weight of
$|\mu|$ over cosets of subgroups of $\SL_2(q_1)$ that have proper
projection in each divisor of $q_1$.  Given $\kappa>0$, there is
$\kappa'>0$ such that if $\mu$ satisfies  for all $q_1 | q$ \beq
\label{e210q} |||\pi_{q_1}(\mu)|||_{\infty} < q_1^{-\kappa}
\|\mu\|_{1} \eeq then \beq \label{e211q} \|\mu \ast \varphi\|_2
\leq C q^{-\kappa'} \|\mu\|_{1} \|\varphi\|_2 \quad \text{for}
\quad \varphi \in E_q. \eeq
\end{lemma}

\section{Bounds for congruence transfer operator}\label{s3}
Our goal is to obtain a bound for powers of transfer operator
 $\| \mM_z^m \|_{\rho}$ for  the family of congruence
subgroups.  Recall that $\mM_z$ acts on functions on $\Sigma
\times G$, so in order to apply Lemma \ref{l1q} we need to
decouple the variables.  Returning to \eqref{e12}, fix $m \leq r <
n$ such that $m=n-r \sim \log q$ to be specified. Write
 $$\mM_{z}^n f(x,g)=$$
\begin{gather}\sum_{i_1, \dots, i_n =1}^l e^{z (\gt(i_n,\dots , i_1,  x)
+ \gt(i_{n-1},\dots , i_1, x) + \dots + \gt(i_1, x))} f((i_n \dots
i_{n-r+1},  0); g_{i_n} \dots g_{i_1} g)\label{e31}\\ +
O(\lambda_{\Re z}^n |f|_{\rho} \rho^r), \label{e31'}
\end{gather}  where the error term refers to the
$L_{l^2(G)}^{\infty}(\Sigma)$-norm.

Fix then the matrices $i_n, \dots, i_{n-r+1}$ and consider the
function $\varphi$ on $G$ defined by
$$\varphi(g)= f((i_n, \dots, i_{n-r+1}, 0), g_{i_n} \dots
g_{i_{n-r+1}} g).$$ We assume $f(x, \cdot) \in E_q$ for each $x$;
hence   $\varphi \in E_q$.

Our aim is to apply  Lemma \ref{l1q}  with \beq \label{e32}
\mu=\sum_{i_1, \dots, i_{n-r}} e^{z ( \gt(i_n,\dots , i_1,  x) +
  \dots + \gt(i_1, x))} \delta_{g_{i_{n-r}} \dots g_{i_1}}. \eeq
Thus by \eqref{e14} we have \beq \label{e33} \|\mu\| \lesssim
\lambda_{\Re z}^m e^{ \Re z ( \gt(i_n,\dots , i_{n-r+1},  0) +
  \dots + \gt(i_{n-r+1}, 0))} e^{\Re z
  \frac{|\tau|_{\rho}}{1-\rho}},
 \eeq
where we used the inequality \beq \label{taucond} | \gt(i_n,\dots
, i_{n-r+1}, i_{n-r}, \dots, i_1,  x)  - \gt(i_n,\dots ,
i_{n-r+1}, \dots)| \leq \tau_{\rho} |\rho|^r. \eeq

 We bound $\|\mu\|_{\infty}$,  which amounts to
 estimating
 \beq \label{e34} \frac{1}{\eqref{e33}} \sum_{g_{i_m} \dots g_{i_1} =g}e^{ \Re z ( \gt(i_n,\dots , i_{1},  x) +
  \dots + \gt(i_1, x))}, \eeq
 where $g$ is fixed.  (We use here the fact that the relation $g_{i_m}
 \dots g_{i_1} = g \, \nod q$ is equivalent to $g_{i_m}
 \dots g_{i_1} = g$ because of the restriction on $m$).  Also
 because  of the index restriction on the transition matrix $A$,
 the condition $g_{i_m} \dots g_{i_1} = g$ specifies $(i_m, \dots,
 i_1) \in \Sigma$ so that estimating \eqref{e34} amounts to
 bounding
 \beq \label{e35} \begin{split}
 &\frac{1}{\eqref{e33}}e^{ \Re z ( \gt(i_n,\dots , i_1,  x) +
  \dots + \gt(i_1, x))}\\ &\sim \lambda^{-m}
e^{ \Re z ( \gt(i_m,\dots , i_1,  x) +
 \dots + \gt(i_1, x))} \end{split}
 \eeq
 for fixed $(i_n, \dots i_1, x) \in \Sigma$.  Thus
 $$ \eqref{e35} = \lambda^{-m} \mM_{\Re z}^m \delta_{(i_m, \dots,
 i_1, x)}(x) \leq \lambda^{-m} \|\mM_{\Re z}^m \phi\|_{\infty} $$
 with $\phi$ any function on $\Sigma$ satisfying $\phi(i_m, \dots,
 i_1, x) =1$.

 By Ruelle's Theorem   (Thm. \ref{tr}),
 $$\| \lambda^{-m} \mM^m \phi - (\int \phi d \nu) h\|_{\rho}
 \lesssim (1-\varepsilon_1)^m \|\phi\|_{\rho},$$
 implying
 \beq \label{e36}\|\lambda^{-m} \mM^m \phi\|_{\infty} \leq c
 \left(\int \phi d \nu + (1- \varepsilon_1)^m
 \|\phi\|_{\rho}\right). \eeq

We may now choose $\phi$ suitably, so  as to obtain an estimate
\beq \label{e37} \lambda^{-m} \|\mM^ m \phi\|_{\infty} \lesssim
e^{-cm}. \eeq

Hence, \beq \label{e38} \eqref{e34} , \eqref{e35} \lesssim
q^{-\kappa} .\eeq

More generally, we also need to evaluate
$|||\pi_{q_1}(\mu)|||_{\infty}$ for $q_1 | q$.  It turns out that
the issue reduces to the previous one, using the following
observation (cf \cite{bg}).  Let $H < G$ and $\pi_p(H) < \SL_2(p)$
proper for each $p|q_1$.  Then we can assume the second commutator
of $\pi_p(H)$ to be trivial if $p|q_1$ and hence the second
commutator of $H$ to be trivial $(\nod \, q_1)$.  Take $m_1  < m$,
$m_1 \sim \log q_1$ so as to ensure that words of length $2m_1$
have norm less than $q_1$.  Using properties of the free group
(see \cite{bg}), it follows from the preceding that the number of
$(i_{m_1}, \dots, i_1) \in \Sigma$ such that $g_{i_{m_1}} \dots
g_{i_1} \in a H$ is bounded by $O(m_1^C)$ for some constant $C$.
Hence we may invoke the estimate on \eqref{e35} with $m$ replaced
by $m_1$ to obtain also
$$ |||\pi_{q_1}(\mu)|||_{\infty} < q_1^{-\kappa} \cdot
\eqref{e33}.$$

Applying Lemma \ref{l1q}, it follows that
$$\|\mu \ast \varphi\|_2 \leq q^{-\kappa'} \|\mu\|_1 \|
\varphi\|_2 \leq q^{-\kappa'} \lambda_{\Re z}^m e^{ \Re z (
\gt(i_n,\dots , i_{n-r+1},  0) +
  \dots + \gt(i_{n-r+1}, 0))} e^{\Re z \frac{|\tau|_{\rho}}{1-\rho}} \|f\|_{\infty} $$ or \beq
\label{e39} \begin{split} &\|\sum_{i_1,\dots, i_{n-r}} e^{z (
\gt(i_n,\dots , i_1,  x) +
  \dots + \gt(i_1, x))} f(i_n, \dots, i_{n-r+1}, 0; g_{i_n} \dots g_{i_1}
g) \|_{l^2(G)} \leq\\ & q^{-\kappa'} \lambda_{\Re z}^m e^{ \Re z (
\gt(i_n,\dots , i_{n-r+1},  0) +
  \dots + \gt(i_{n-r+1}, 0))} e^{\Re z \frac{|\tau|_{\rho}}{1-\rho}} \|f\|_{\infty}. \end{split} \eeq

Summing \eqref{e39} over $i_n, \dots, i_{n-r+1}$ implies by
\eqref{e14} again

 \beq \label{e310}\begin{split} & \|
\sum_{i_1, \dots, i_n =1}^l e^{z (\gt(i_n,\dots , i_1,  x) +
\gt(i_{n-1},\dots , i_1, x) + \dots + \gt(i_1, x))} f((i_n \dots
i_{n-r+1},  0); g_{i_n} \dots g_{i_1} g)\|_{l^2(G)}\\ &\lesssim
q^{-\kappa'}  \lambda^{m+r} \|f \|_{\infty}. \end{split} \eeq

Therefore it follows that if $n > \log q$ \beq \label{e311}
\begin{split} \|\mM_z^n f\|_{L^{\infty}_{l^2(G)}(\Sigma)} & \leq
\lambda_{\Re z}^n (q^{-\kappa'} \|f\|_{\infty} + \rho^r
|f|_{\rho}) \\& \leq \lambda_{\Re z}^n q^{-\kappa'}
(\|f\|_{\infty} + \rho^{\frac{n}{2}} |f|_{\rho})\end{split} \eeq
for \beq \label{e312} f \in \mathcal{F}'_{\rho} =
\mathcal{F}_{\rho} \cap \mathcal{C}_{E_q}(\Sigma). \eeq Note that
in \eqref{e311} there is no restriction on $\Im z$.

We also need to estimate $|\mM_{z}^n f|_{\rho}$.

Let $x, y \in \Sigma$ be such that $x_i = y_i$ for $0 \leq i < l$.
Estimate  $$ |\mM_{z}^n f(x, g) - \mM_{z}^n f(y, g)|\leq$$
\begin{gather}
\sum_{i_1, \dots, i_n} e^{\Re z (\tau(i_n, \dots, i_1, x) + \dots
+ \tau(i_1, x))}| f(i_n, \dots, i_1 x; g_{i_n} \dots g_{i_1} g) -
f(i_n, \dots, i_1 y; g_{i_n} \dots g_{i_1} g)| \label{e313} \\  +
|\sum_{i_1, \dots, i_n} \left( e^{z (\tau(i_n, \dots, i_1, x) +
\dots + \tau(i_1, x))} -e^{z (\tau(i_n, \dots, i_1, y) + \dots +
\tau(i_1, y))}\right) f(i_n, \dots, i_1 y; g_{i_n} \dots g_{i_1}
g)|. \label{e314}
\end{gather}

Clearly for the first term we have \beq \label{e315} \eqref{e313}
\lesssim \lambda_{\Re z}^n |f|_{\rho} \rho^{n+l}. \eeq

To estimate \eqref{e314} we repeat the argument leading to
\eqref{e311}.  Thus we bound \eqref{e314} as follows \beq
\label{e316} |\sum_{i_1, \dots, i_n} \left( e^{z (\tau(i_n, \dots,
i_1, x) + \dots + \tau(i_1, x))} -e^{z (\tau(i_n, \dots, i_1, y) +
\dots + \tau(i_1, y))}\right)f(i_n, \dots, i_{n-r+1}, 0; g_{i_n}
\dots g_{i_1} g)| \eeq \beq \label{e317} + \rho^{r} |f|_{\rho}
\sum_{i_1, \dots, i_n} |\left( e^{z (\tau(i_n, \dots, i_1, x) +
\dots + \tau(i_1, x))} -e^{z (\tau(i_n, \dots, i_1, y) + \dots +
\tau(i_1, y))}\right)|. \eeq

Estimate \beq \label{e318} \begin{split} & \sum_{i_1, \dots, i_n}
|\left( e^{z (\tau(i_n, \dots, i_1, x) + \dots + \tau(i_1, x))}
-e^{z (\tau(i_n, \dots, i_1, y) + \dots + \tau(i_1, y))}\right)|
\\& \leq \sum_{i_1, \dots, i_n} e^{\Re z (\tau(i_n, \dots, i_1, x)
+ \tau(i_1, x))} |1 - e^{z (\tau(i_n, \dots, i_1, y) + \dots +
\tau(i_1, y)- \tau(i_n, \dots, i_1, x) -\dots - \tau(i_1, x) )}|
\\& \lesssim \lambda_{\Re z}^n (1+ | \Im z|) |\tau|_{\rho}
(\rho^{n+l} + \dots + \rho^{1+l}) < \lambda_{\Re z}^n \frac{1 +
|\Im z|}{1 -\rho} |\tau|_{\rho} \rho^l;
\end{split} \eeq
therefore \beq \label{e319} \eqref{e317} < \lambda^n \frac{1+ |\Im
z |}{1-\rho} |\tau|_{\rho} \rho^{l+r} |f|_{\rho}. \eeq

To bound \eqref{e316} we apply again the convolution estimate on
$G$ from section \ref{s2}.  Consider the measure \beq \label{e320}
\nu = \sum_{i_1, \dots, i_{n-r}} \left( e^{z (\tau(i_n, \dots,
i_1, x) + \dots + \tau(i_1, x))} -e^{z (\tau(i_n, \dots, i_1, y) +
\dots + \tau(i_1, y))}\right)\delta_{g_{i_{n-r}} \dots g_{i_1}}
\eeq with $i_n, \dots i_{n-r+1}$ fixed.

Repeating \eqref{e318} gives (with $m=n-r$) \beq \label{e321}
\|\nu\| \lesssim \lambda^m (1+ | \Im z |) |\tau|_{\rho}
\frac{\rho^l}{1-\rho} e^{\Re z (\tau(i_n, \dots, i_{n-r+1}, 0)
+\dots + \tau(i_{n-r+1}, 0))}. \eeq

Also, as above, we have \beq \label{e322} \begin{split}
&\frac{1}{\eqref{e321}} \|\nu\|_{\infty} =
\frac{1}{\eqref{e321}}|e^{z(\tau(i_n, \dots, i_1, x)+ \dots +
\tau(i_1, x))}- e^{z(\tau(i_n, \dots, i_1, y)+ \dots + \tau(i_1,
y))}|\\&\leq \lambda^{-m} e^{\Re z(\tau(i_n, \dots, i_1, x)+ \dots
+ \tau(i_1, x))} \overset{\eqref{e38}}{<} q^{-\kappa}.
\end{split} \eeq
and \beq \label{e322a} \frac{1}{\eqref{e321}} |||
\pi_{q_1}(\nu)|||_{\infty} < q_1^{-\kappa} \quad \text{for} \quad
q_1|q. \eeq

Therefore, by the results from section \ref{s2}, we obtain (with
$\varphi$ defined as in section \ref{s2}) that \beq
\label{e323}\begin{split}& \| \nu \ast \varphi \|_{l^2(G)} \leq
q^{-\kappa'}  \eqref{e321} \|f\|_{\infty} \lesssim\\
&q^{-\kappa'}\lambda^m (1+ | \Im z |) |\tau|_{\rho}
\frac{\rho^l}{1-\rho} e^{\Re z (\tau(i_n, \dots, i_{n-r+1}, 0)
+\dots + \tau(i_{n-r+1}, 0))} \|f\|_{\infty}.
\end{split}  \eeq Summation over $i_n, \dots, i_{n-r+1}$ gives then

\beq \label{e324} \|\eqref{e316} \|_{l^2(G)} \lesssim q^{-\kappa'}
\lambda^n (1+| \Im z |) | \tau|_{\rho} \rho^l \| f\|_{\infty}.
\eeq

From \eqref{e315}, \eqref{e319}, \eqref{e324}, it follows that
\beq \label{e325} \| \mM_z^n f(x, \cdot) -\mM_z^n f(y,
\cdot)\|_{l^2(G)} \lesssim \rho^l \lambda_{\Re z}^n \left( \rho^n
|f|_{\rho} + \rho^r (1 + | \Im z|) |f|_{\rho} + q^{-\kappa'} (1 +
| \Im z|) \|f\|_{\infty} \right). \eeq

Therefore, if $n > \log q$ we have \beq \label{e326} |\mM_z^n
f|_{\rho} \leq C \lambda_{\Re z}^n q^{-\kappa'} \left(
\|f\|_{\infty} + \rho^{n/2} |f|_{\rho} \right) (1+ | \Im z|).\eeq

Take $n$ such that \beq \label{e328} n \sim \log q + C \log (1+ |
\Im z|) \eeq for a suitable constant $C$.  It follows from
\eqref{e311}, \eqref{e326} that \beq \label{e329} \| \mM_z^n f
\|_{\infty} + \rho^{\frac{n}{2}} | \mM_z^n f|_{\rho} <
\lambda_{\Re z}^n q^{-\kappa'} ( \|f\|_{\infty} +
\rho^{\frac{n}{2}} |f|_{\rho}). \eeq

Iterating \eqref{e329} shows that if $f \in \mathcal{F}'_{\rho}$,
then for all $m \in \zed_+$ \beq \label{e330} \| \mM_z^{mn}
f\|_{\infty} + \rho^{\frac{n}{2}}|\mM_{z}^{mn} f|_{\rho} \leq
\lambda_{\Re z}^{mn} q^{-m \kappa'} \|f\|_{\rho}, \eeq and hence
\beq \label{e331} \|\mM_z^{mn} f \|_{\rho} < \lambda_{\Re z}^{mn}
q^{-m \kappa'}  q (1+|\Im z|) \|f\|_{\rho}, \eeq where $n$ is
given by \eqref{e328}.  Thus for $m \geq 1$ \beq \label{e332} \|
\mM_z^m |_{\mathcal{F}'_{\rho}} \|_{\rho} < \lambda_{\Re z}^m q^{-
\frac{m}{n} \kappa'} q (1+|\Im z|). \eeq We distinguish two cases:
$\log (1 + | \Im z|) \lesssim \log q$ and $\log q  \ll \log (1 + |
\Im z|)$. The conclusion is the following:
\begin{lemma} \label{props3} Notation being as above, there is
$\varepsilon > 0$ such that \beq \label{e333} \| \mM_z^m
|_{\mathcal{F}'_{\rho}} \|_{\rho}  < q^C e^{-\varepsilon m}
\lambda_{\Re z}^m \, \, \text{if} \, \, |\Im z| \leq q \eeq and
\beq \label{e334} \| \mM_z^m |_{\mathcal{F}'_{\rho}} \|_{\rho}  <
|\Im z|^C e^{-\varepsilon \frac{\log q}{\log | \Im z|} m}
\lambda_{\Re z}^m \, \, \text{if} \, \, |\Im z| > q. \eeq
\end{lemma}

\section{Resolvent of congruence transfer operator} \label{s6}

We now use Lemma \ref{props3}  to estimate the resolvent
$(I-\mM_z)^{-1}$ on $\mathcal{F}'_{\rho}$. By \eqref{e6.6},
\eqref{e6.7}  this will provide us bounds on $F(z,x)$, assuming
$\varphi \in E_q$.

 Take $\Re z < -\delta + \varepsilon_1$ such
that \beq \label{e61} \lambda_{-\delta+ \varepsilon_1} <
e^{\frac{\varepsilon}{2}} \eeq with $\varepsilon > 0$ from
\eqref{e333}.  If $| \Im z| < q$, we obtain \beq \label{e62} \|
(I- \mM_z)^{-1}|_{\mF'_{\rho}}\| < q^C \sum e^{-\vge m}
\lambda_{\Re z}^m \lesssim \frac{1}{\vge}q^C. \eeq

If $| \Im z | \geq q$, we impose the restriction \beq \label{e63}
\Re z < -\gd + \vge_2 \frac{\log q}{\log| \Im z |} \eeq with
$\vge_2 > 0$ small enough to ensure that
$$\lambda_{\Re z} e^{-\vge \frac{\log q}{\log| \Im z
|}} <  e^{-\frac{\vge}{2} \frac{\log q}{\log| \Im z |}}.$$ Under
this restriction on $z$, we obtain from \eqref{e334} that \beq
\label{e64} \| (I- \mM_z)^{-1}|_{\mF'_{\rho}}\| < | \Im z |^C.
\eeq

In summary, we proved the following

\begin{theorem} \label{p1}
The resolvent $(I- \mM_z)^{-1}|_{\mF'_{\rho}}$ is holomorphic on
the complex region $D(q)$ given by \beq \label{e65} \Re z < -
\delta + \vge_2 \min\left(1, \frac{\log q}{\log(|\Im z| +1)}
\right) \eeq (with $\vge_2$ independent of $q$) and satisfies the
estimate \beq \label{e66}
 \| (I- \mM_z)^{-1}|_{\mF'_{\rho}}\| <(q+ | \Im z |)^C.
\eeq
\end{theorem}

Returning to \eqref{e6.6}, \eqref{e6.7}, it follows that for
$\varphi \in E_q$ the Laplace transform $F(z,x)$ of $N(a,x)$ is
bounded by
 \beq \label{e67} | F(z, x)| \lesssim \frac{(q + |\Im z|)^C}{|z|}
\| \varphi\|_2\eeq for $z$ satisfying \eqref{e65}.

To extract information about $N(a)$ we apply Fourier inversion to
\eqref{e58}, following the argument in \cite[p. 31]{l89} (but with
a different class of functions $k$).

Specify some  smooth and compactly supported bump function $k$ on
$\reals$.  We get from \eqref{e58}

\beq \label{e68} \int_{-\infty}^{\infty} k(t) e^{-\delta t} N(a+t)
dt = e^{\delta a} \int e^{-i a \theta} \hat{k}(-i \theta)
F(-\delta+ i \theta) d \theta, \eeq where
$$\hat{k}(z) = \int e^{ z t} k(t) d t $$ is an entire function.

Note that $|\hat{k}(i \theta)|$ is rapidly decaying since $k$ is
smooth.

In fact, proceeding more precisely, fix a small parameter $\gamma
> 0$ (the localization of $k$) and consider functions
\beq \label{e69} k_{\gamma}(t) = \frac{1}{\gamma}
K\left(\frac{t}{\gamma}\right), \eeq where $K$ is a fixed smooth
bump function such that \beq \label{e610} \int K =1, \eeq \beq
\label{e611} \supp K \subset [-\frac{1}{2}, \frac{1}{2}],  \eeq
\beq \label{e612} |\hat{K}(\gl)| \lesssim e^{-|\gl|^{\frac{1}{2}}}
\quad \text{for} \quad |\gl| \to \infty. \eeq Hence \beq
\label{e613} |\hat{k}(z)| \lesssim e^{-|\gamma |z||^{\frac{1}{2}}}
\quad \text{for} \quad |\Re z| < O(1). \eeq

Returning to \eqref{e68}, modify the line of integration $\Re z =
0$ to the curve $$z(\theta) = w(\theta) + i \theta ,$$ where \beq
\label{e614} w(\gt) = \frac{1}{2} \vge_2 \min \left( 1, \frac{\log
 q}{\log( 1 + |\theta|)}\right), \eeq
so as to remain in the analyticity region given by Theorem
 \ref{p1}.

 We obtain
 $$e^{\delta a} \int_{-\infty}^{\infty} e^{-a z(\theta)}
 \hat{k}(z(\theta)) F(-\delta+ z(\theta)) d \theta =
e^{\delta a} \int_{-\infty}^{\infty} e^{-a w(\theta) - i a \theta}
 \hat{k}(w(\theta) + i \theta ) F(-\delta+ w(\theta) + i \theta ) d
 \theta,$$
which is bounded by \beq \label{e615} \|\varphi\|_2 e^{\delta a }
\int_{-\infty}^{\infty} e^{-a w(\theta)} e^{-(\gamma
|\theta|)^{\frac{1}{2}}} (q+ |\theta|)^C d \theta, \eeq applying
\eqref{e613} and \eqref{e67}.

From the definition of $w(\theta)$ it is clear that \beq
\label{e616} \eqref{e68}, \eqref{e615} < e^{\delta a} q^{C}
\gamma^{-C} \exp \left(- a \vge_3 \min \left( 1, \frac{\log
q}{\log \frac{a}{\gamma}}\right) \right) \|\varphi\|_{2}. \eeq

This proves (replacing $k(t)$ by $ e^{\delta
 t} k(t)$)

 \begin{prop} \label{p2}
Let $\varphi \in E_q$ and $N(a)$ given by \eqref{e56}. Then \beq
\label{e617} \left | \int_{-\gamma/2}^{\gamma/2} k_{\gamma}(t)
N(a+t) d t\right| <  q^{C} \gamma^{-C} \exp \left(- a \vge_3 \min
\left( 1, \frac{\log q}{\log \frac{a}{\gamma}}\right) \right)
e^{\delta a} \|\varphi\|_{2}. \eeq
\end{prop}

\section{Bound for the error term} \label{s7}

Next consider the case where in \eqref{e54}, $\phi =1$ (the
constant function).

Here we consider simply the action of $\mL_z$ on
$\mF_{\rho}(\Sigma)$ exactly as in \cite{l89}, but we use the
stronger estimates on $(I-\mL_z)^{-1}$ following from \eqref{e42},
given by \cite{naud05a}.

If $\Re z < -\delta+ \vge_4$ (with $\vge_4$ small enough) and $z
\notin U$  (some complex neighborhood of $-\delta$), \eqref{e42}
implies that \beq \label{e71} \|(I- \mL_z)^{-1}\| \lesssim | \Im
z|^{2}. \eeq

For $s \in U$, apply \eqref{e43}.  Thus
$$ \mL_z^n = \lambda_z^n(\nu_z \otimes h_z) + (\mL''_z)^n,$$
where $\|(\mL''_z)^n\| < e^{\vge n}$ by \eqref{e44}.

Hence for $z \in U$ \beq \label{e72}
 (I-\mL_z)^{-1} = \frac{1}{1-\lambda_z}(\nu_z \otimes h_z) +
 (I-\mL''_z)^{-1} \eeq
 with $(I-\mL''_z)^{-1}$ holomorphic  (this is Proposition 7.2 in
 \cite{l89}).

 Combining \eqref{e71}, \eqref{e72} we get
 \begin{prop} \label{p3}  Consider $\mL_z$ acting on
 $\mF_{\rho}(\Sigma)$.  Then for $\Re z <
 -\delta+ \vge_5$
 \beq \label{e73}
 (I-\mL_z)^{-1} - \frac{1}{1-\lambda_z}(\nu_z \otimes h_z)\eeq is
 holomorphic and bounded by $C ( | \Im z|^2 + 1)$.

 Let $\mu_z= \nu_z \otimes h_z$.  The function $\frac{1}{1- \lambda_z}$ has a pole at $z = -
 \delta$ with residue $$-\frac{1}{(\frac{d}{ d z}
 \lambda_z)|_{z=-\delta}} = \frac{1}{\int \tau d \mu_{-\delta}}.
 $$
 Consequently
 \beq \label{e74}
 (I-\mL_z)^{-1} - \frac{\nu_{-\delta} \otimes h_{-\delta}}{\int \tau d \mu_{-\delta}} \frac{1}{z+\delta} \eeq
is analytic for $\Re z < - \delta + \vge_5$.
\end{prop}

Letting \beq \label{e75} N(a)=\sum_{\substack{ y \in \Lambda \\
d_H(0, y w)-d_H(0, w) \leq a}} 1 \eeq and \beq \label{e76}
F(z)=\int e^{az} N(a) d a, \eeq it follows from \eqref{e510} ,
\eqref{e74} that \beq \label{e77} F(z) = \frac{1}{z}
(I-\mL_z)^{-1} 1 = \frac{h_{-\delta} (\xi \equiv w)}{\int \tau d
\mu_{-\delta}} \frac{1}{z+ \delta} + G(z),\eeq where $G(z)$ is
analytic on $\Re z < -\delta + \vge_5$ and bounded by $C (| \Im
z|^2 +1)$.

As in section \ref{s6} we have \beq \label{e78} \begin{split} &
\int k_{\gamma}(t) e^{-\delta t} N(a+t) d t  = e^{\delta a}\int
\hat{k}_{\gamma}(i \theta) F(-\delta + i \theta) e^{- i a \theta} d \theta \\
& =e^{\delta a} \left( C_0 \int_{PV}e^{-ia \theta}
\hat{k}_{\gamma}( i \theta) \frac{1}{i \theta} d \theta + \int
e^{- i a \theta} \hat{k}_{\gamma} (i \theta) G(-\delta + i \theta)
\right),
\end{split} \eeq
where \beq \label{e79} C_0 = \frac{h_{-\delta}(\xi \equiv w)}{\int
\tau d \nu_{- \delta}}. \eeq

The second term in \eqref{e78} is estimated by moving the line of
integration $\Re z = 0$ to $\Re z = \frac{1}{2} \vge_5$.  We
obtain by \eqref{e613} and the assumption on $G$: \beq
\label{e710} |\int e^{- i a \theta} \hat{k}_{\gamma} (i \theta)
G(-\delta + i \theta)|\lesssim e^{-\frac{\vge_5}{2} a} \int (1+
\theta^2) e^{- (\gamma |\theta|)^{\frac{1}{2}}} d \theta < c
\gamma^{-3} e^{-\frac{\vge_5}{2} a}. \eeq Also
$$\int_{PV}e^{-ia \theta} \hat{k}_{\gamma}( i \theta)
\frac{1}{i \theta} d \theta = \int_{0}^{\infty} k_{\gamma}(t+a) d
t \overset{\eqref{e610}}{=} 1.$$

Therefore we obtain
\begin{prop} \label{p4} Let $N(a)$ be given by \eqref{e75}.Then
\beq \label{e711} \int k_{\gamma}(t) N(a+t) dt = C_0 e^{\delta a}
+ o\left(\gamma^{-3} e^{(\delta - \vge_6) a} \right) \eeq for some
$\vge_6>0$.  Here $C_0$ is a fixed constant.
\end{prop}

\section{Proof of Theorem \ref{t5a}}  \label{s7a}
Let $\varphi$ be a function on $\SL_2(q)$ and let $N(a,x)$  denote
the counting function given as above by
$$ N(a,x) = \sum_{y \in \Lambda, d_H(0, y  x w) - d_{H}(0, x w)
\leq a} \varphi(\pi_q(y)).$$ What we proved in Theorem \ref{p1} is
that for $\varphi \in E_q$ the Laplace transform of $N(a,x)$ in
$a$  (given by \eqref{e58}) is holomorphic on $D(q)$ given by \beq
\label{e65aa} D(q) = \left\{z \, : \, \Re z < - \delta + \vge_2
\min\left(1, \frac{\log q}{\log(|\Im z| +1)} \right) \right \}
\eeq with $\vge_2$ independent of $q$. Let us denote by $\mL_z(q)$
the dynamical transfer operator on the congruence subgroup
$\Lambda(q)$.  Thus $\det(1-\mL_z(q))$ is the dynamical (Ruelle's)
zeta-function  associated with  the congruence subgroup
$\Lambda(q)$.   Using \eqref{e510} we have that the Laplace
transform of $N(a,x)$ is also obtained as the inverse of
$(1-\mL_z(q))$.   Now considering the action of $\mL_z(q)$ on
$\mF_{\rho}(\Sigma(q))$, recalling the decomposition of
$L^2(\SL_2(q))$ given by \eqref{e82}, and applying theorem
\ref{p1} to $E_{q_1}$ for all $q_1 |q$, and Proposition \ref{p3}
to the constant function, we obtain that $1-\mL_z(q)$ has
holomorphic inverse (apart from $z = -\delta$) on
$$D= D(1) \cap \bigcap_{q_1 | q} D(q_1),
$$ where $$D(1) = \{z \, : \, \Re z < -\delta + \vge_5 \}$$ by
Proposition \ref{p3}.  Consequently $D$ is given by \beq
\label{e65cor} D = \left\{ z\, : \, \Re z < - \delta + \vge_6
\min\left(1, \frac{1}{\log(|\Im z| +1)} \right) \right \}\eeq for
some $\vge_6$ independent of $q$,
 implying that dynamical  zeta function $\det(1-\mL_z(q))$
 has no zeros on $D$ (apart from simple zero at
$-\delta$).

Theorem \ref{t5a} now follows from the  equality of the dynamical
zeta-function  and Selberg zeta function (Theorem 15.8 in
\cite{bort}), and the correspondence between the zeros of
Selberg's zeta function and resonances (see  \cite{patper} and
Chapter 10 in \cite{bort}).

\section{Proof of Theorem \ref{t5b}} \label{s8}

Propositions \ref{p3} and \ref{p4} are our basic estimates used in
the proof of Theorem \ref{t5b}.  Note that what comes out
Proposition \ref{p2} will only play the role of error terms.

Fix a modulus $q$, $(q, q_0) =1$ (with $q_0$ given by the strong
approximation property) and $q$ square-free.

For some element $\xi \in \SL_2(q)$ we need to evaluate \beq
\label{e81} N(a; q, \xi)= |\{ y \in \Lambda; \pi_q(y)= \xi \quad
\text{and} \quad d_{H}(0, y w) \leq a\}| \eeq(replace $a$ by
$a+d_H(0, w)$ in \eqref{e56}).

Recall the decomposition of the space $L^2(\SL_2(q))$ in
\eqref{e82}.  Writing \beq \label{e83} 1_{g=\xi}=
\frac{1}{|\SL_2(q)|}+ \sum_{\substack{q_1|q \\ q_1 \neq
1}}P_{q_1}(1_{g=\xi}), \eeq we get
\begin{gather}N(a; q, \xi)= \frac{1}{|\SL_2(q)|}\sum_{\substack{y \in
\Lambda \\ d_H(0, y w) \leq a}} 1 \label{e84} \\ +
\sum_{\substack{q_1|q
\\ q_1 \neq 1}} \sum_{\substack{y \in \Lambda \\ d_H(0, y w) \leq
a}} \varphi_{q_1}(\pi_{q_1}(y)) \label{e85}
\end{gather}
with $$\varphi_{q_1}=P_{q_1}(1_{g=\xi}) \in E_{q_1}.$$

Thus \beq \label{e86} \|\varphi_{q_1}\|_{2} <
\frac{|\SL_2(q_1)|^{\frac{1}{2}}}{|\SL_2(q)|}.\eeq

We use Proposition  \ref{p4} to evaluate the right-hand side of
\eqref{e84} and Proposition \ref{p2} to bound terms \eqref{e85}.
Hence, fixing some $\gamma > 0$
\begin{gather}
\int k_{\gamma}(t) N(a+t; q, \xi) d t = \frac{C_1}{|\SL_2(q)|}
e^{\delta a} + o(\gamma^{-3} e^{-c a} e^{\delta a} )\label{e87}\\+
\gamma^{-C} \sum_{\substack{q_1|q
\\ q_1 \neq 1}} q_1^C \exp\left(- c a \min\left(1, \frac{\log q_1}{\log
\frac{a}{\gamma}}
\right)\right)\frac{|\SL_2(q_1)|^{\frac{1}{2}}}{|\SL_2(q)|}e^{\delta
a}. \label{e88}
\end{gather}

We estimate \eqref{e88} as \begin{gather} \frac{\gamma^{-C}
|a|^C}{|\SL_2(q)|}\left( \sum_{\substack{q_1 | q \\ 1 < q_1 <
\frac{|a|}{\gamma}}}e^{- c a \frac{\log q_1}{\log
\frac{a}{\gamma}}}\right)e^{\delta a} \label{e89} \\+ \gamma^{-C}
q^{C} e^{(\delta-c) a} \label{e810} \end{gather} and \beq
\label{e811} \begin{split}  & \sum_{\substack{q_1 | q \\ q_1 \neq
1 }}e^{- c a \frac{\log q_1}{\log \frac{a}{\gamma}}} < \prod_{p|q}
\left(1 +e^{- c a \frac{\log p}{\log \frac{a}{\gamma}}}\right) - 1 \\
&< \exp\left( \sum_{s=2}^{\infty} e^{-ca \frac{\log s}{\log
\frac{a}{\gamma}}} \right)-1 < e^{-c \frac{a}{\log
\frac{a}{\gamma}}}, \end{split}\eeq assuming\beq \label{e812} \log
\frac{1}{\gamma} \ll a. \eeq Therefore we proved the following
result, of which Theorem \ref{t5b} is an immediate consequence.
\begin{prop} \label{pt1}Notation being as above, we have
\beq \label{e813} \int k_{\gamma}(t) N(a+t; q, \xi) d t=
\frac{e^{\delta a}}{|\SL_2(q)|} \left(C_1 + o(\gamma^{-C} e^{-c
\frac{a}{\log \frac{a}{\gamma}}})\right) + \gamma^{-C} q^{C}
e^{(\delta-c) a}, \eeq where we assume \beq \label{e814} \log
\frac{1}{\gamma} \ll \frac{a}{\log a}. \eeq
\end{prop}

We remark that what is really required for sieving applications is
a bound for the ratio \beq \label{ratio} \frac{\int k_{\gamma}(t)
N(a+t; q) d t}{\int k_{\gamma}(t) N(a+t) d t} \eeq of the form
$$\frac{1}{|\SL_2(q)|} + O(e^{-\varepsilon a} q^C)$$ or
$$\frac{1}{|\SL_2(q)|}\left( 1+ O(e^{-c \frac{a}{\log a}})\right) + O(e^{-\varepsilon a} q^C).$$
To bound the ratio \eqref{ratio} it suffices to use Proposition
\ref{p2} (which builds crucially on the generalized expansion
result given by Lemma \ref{l1q}) combined with the result of
Lalley \cite{l89}. Of course, the results of Dolgopyat \cite{dolg}
and Naud \cite{naud05a}, \cite{naud05b} are necessary to establish
Theorem \ref{t5a}, which is of independent interest.

\section{Proof of Theorem \ref{t2}} \label{sec:t2}

\subsection{Combinatorial sieve} As in \cite{BGS}, we will make use of the simplest
combinatorial sieve which is turn is based on the Fundamental
Lemma in the theory of elementary sieve, see \cite{iw} and
\cite{hr}. Our formulation is tailored for the applications below.

Let $A$ denote a finite sequence $a_n$, $n\ge 1$ of nonnegative
numbers. Denote by $X$ the sum \beq \label{cs1} \sum_{n}a_n=X.
\eeq $X$ will be large, in fact tending to infinity. For a fixed
finite set of primes $B$ let $z$ be a large parameter (in our
applications $z$ will be a small power of $X$ and $B$ will usually
be empty). Let \beq \label{cs2} P=P_z=\prod_{\substack{p\le z \\p
\notin B}}p . \eeq Under suitable assumptions about sums of $A$
over $n$'s in progressions with moderate-size moduli $d$, the
sieve gives upper and lower estimates which are of the same order
of magnitude for sums of $A$ over the $n$'s which remain after
sifting out numbers with prime factors in $P$.

More precisely, let \beq \label{cs3} S(A, P) := \sum_{(n, P)=1}
a_n. \eeq The assumptions on sums in progressions are as follows:

($A_0$) For $d$ square-free, and having no prime factors in $B$
($d < X$), we assume that the sums over multiples of $d$ take the
form \beq \label{cs4} \sum_{n \equiv 0(d)} a_n = \beta(d)X +r(A,
d), \eeq where $\beta(d)$ is a multiplicative function of $d$ and
$$ \text{for}\, p \notin B,  \beta(d) \leq 1-\frac{1}{c_1} \, \text{for a
fixed}\, c_1.$$ The understanding being that $\beta(d) X$ is the
main term and that the remainder $r(A, d)$ is smaller, at least on
average (see the next axiom).

\medskip

 ($A_1$) $A$ has level distribution
$D =D(X)$, ($D < X)$ that is
$$\sum_{d \le D}|r(d, A)| \ll \frac{X}{(\log X)^B} \, \, \text{for
all} \, \, B>0.
$$

\medskip

($A_2$) $A$ has sieve dimension $t>0$, that is for a fixed $c_2$
we have
$$ \left|\sum_{\substack{w \le p \le z\\p \notin B}}\beta(p) \log p - t
\log \frac{z}{w} \right| \le c_2$$ for $2 \le w \le z$.

\medskip

In terms of these conditions ($A_0$),  ($A_1$), ($A_2$) the
elementary combinatorial sieve yields:
\begin{theorem}\label{t:cs} Assume ($A_0$),
($A_1$) and ($A_2$) for  $s > 9 t$ and $z=D^{1/s}$ and $X$ large
we have \beq \label{cs} \frac{X}{(\log X)^t} \ll S(A, P_z) \ll
\frac{X}{(\log X)^t}. \eeq
\end{theorem}

\subsection{Applying the sieve} Now let $\Lambda$ be a Zariski dense subgroup of  $\SL(2,
\zed)$ and let $f \in \zed[x_{ij}]$ be weakly primitive with
$t(f)$ irreducible factors. The key nonnegative sequence $a_n$ to
which we apply the combinatorial sieve is defined as follows: for
$n \geq 0$ we let  \beq
\label{anseq}  a_n=a_n(T)= \sum_{\substack{\gamma \in \Lambda \\
|\gamma| \leq T\\|f(\gamma)| = n}}1. \eeq

The sums on progressions are then, for $d \geq 1$ square free \beq
\label{consum1} \sum_{n \equiv 0 (d)} a_n(T) =
\sum_{\substack{\gamma \in \Lambda; |\gamma| \leq T\\ f(\gamma)
\equiv 0 (d)}}1 = \sum_{\substack{\rho \in \Lambda / \Lambda(d) \\
f(\rho) \equiv 0 (d)}} \sum_{\substack{\gamma \in \Lambda(d) \\
|\gamma| \leq T}}1. \eeq Consider the case of $\delta < 1/2$; the
case of $\delta > 1/2$ is similar and simpler. According to
Theorem \ref{t5b}  we then obtain \beq \begin{split} &\sum_{n
\equiv 0 (d)} a_n(T) =  \sum_{\substack{\rho \in \Lambda / \Lambda(d) \\
f(\rho) \equiv 0 (d)}} \frac{T^{ 2 \delta}}{| \Lambda_d|}\left(1 +
O\left(T^{-\frac{1}{\log \log T}}\right)\right) + O \left(d^C T^{2
\delta - \varepsilon_1}\right)\\&= X
\frac{|\Lambda_d^f|}{|\Lambda_d|} +
O\left(\frac{|\Lambda_d^f|}{|\Lambda_d|} X^{1-\frac{1}{ 2 \delta
\log \log X}}\right) + O\left( |\Lambda_d^f| d^C X^{1-
\frac{\varepsilon_1}{2\delta}}\right).\end{split}, \eeq where \beq
X= \sum_{k
\in \natls} a_k(T) = \sum_{\substack{\gamma \in \Lambda \\
|\gamma| \leq T}}1, \eeq  $\Lambda_d$ is the reduction of
$\Lambda$ mod $d$, and $\Lambda_d^f$ is the subset of $\Lambda_d$
at which $f(x) = 0 \, \nod \, d$.

Using strong approximation theorem \cite{mvw} and Goursat lemma as
in \cite{BGS} we obtain  that outside of finite set of primes
$S(\Lambda)$  we have $\Lambda_p  \cong \SL_2(\fp)$ and $\Lambda
\to \Lambda_{d_1} \times \Lambda_{d_2}$ is surjective for $(d_1,
d_2) =1$ and $d_1 d_2$ square free.  Let \beq
\beta(d)=\frac{|\Lambda_d^f|}{|\Lambda_d|}. \eeq Using Lang-Weil
theorem \cite{lw} as in \cite{BGS}, we obtain \beq \label{ldfb}
|\Lambda_d^f| \ll d^2 \eeq and \beq \label{ontr}
\frac{|\Lambda_d^f|}{|\Lambda_d|} = \frac{t(f)}{p} +
O(p^{-\frac{3}{2}}). \eeq

 Hence we have \beq \sum_{n \equiv 0
(d)} a_n(T) = \beta(d) X + r(A,d), \eeq with \beq r(A,d) \ll
\frac{1}{d} X^{1-\frac{1}{ 2 \delta \log \log X}} + d^{C+2}X^{1-
\frac{\varepsilon_1}{2\delta}}.\eeq  Verification of $(A_0)$ is
completely analogous to Proposition 3.1 in \cite{BGS}.  Regarding
the level distribution $(A_1)$  we have that \beq \sum_{d \leq D}
|r(A,d)| \ll X^{1-\frac{1}{ 2 \delta \log \log X}} + D^{C+3}X^{1-
\frac{\varepsilon_1}{2\delta}} \ll \frac{X}{(\log X)^B} \eeq for
any $B>0$ as long as \beq D \leq X^\tau \, \, \text{with} \, \,
\tau < \frac{2 \delta}{(C+3) \varepsilon_1}. \eeq

Finally, to verify the third axiom concerning the sieve dimension,
we have, using \eqref{ontr}, that   \beq  \sum_{w \le p \le z}
\beta(p) \log p = \sum_{w \le p \le z}\left( \frac{t \log p}{p} +
O\left(\frac{\log p}{p^{3/2}}\right)\right) = t \log \frac{z}{w}
+O(1), \eeq which establishes $(A_2)$ with the sieve dimension
being $t$.

\end{document}